# ESTIMATION IN SEMIPARAMETRIC SPATIAL REGRESSION[1]

BY JITI GAO, ZUDI LU[2] AND DAG TJØSTHEIM

*University of Western Australia, Chinese Academy of Sciences
and University of Bergen*

Nonparametric methods have been very popular in the last couple of decades in time series and regression, but no such development has taken place for spatial models. A rather obvious reason for this is the curse of dimensionality. For spatial data on a grid evaluating the conditional mean given its closest neighbors requires a four-dimensional nonparametric regression. In this paper a semiparametric spatial regression approach is proposed to avoid this problem. An estimation procedure based on combining the so-called marginal integration technique with local linear kernel estimation is developed in the semiparametric spatial regression setting. Asymptotic distributions are established under some mild conditions. The same convergence rates as in the one-dimensional regression case are established. An application of the methodology to the classical Mercer and Hall wheat data set is given and indicates that one directional component appears to be nonlinear, which has gone unnoticed in earlier analyses.

**1. Introduction.** Data collected at spatial sites occur in many scientific disciplines, such as econometrics, environmental science, epidemiology, image analysis and oceanography. Often the sites are irregularly positioned, but, with the increasing use of computer technology, data on a regular grid and measured on a continuous scale are becoming more and more common. This is the kind of data that we will be considering in this paper.

In the statistical analysis of such data, almost exclusively, the emphasis has been on parametric modeling. So-called joint models were introduced in the papers by Whittle [36, 37], but, after the ground breaking paper by Besag

Received May 2003; revised July 2005.
[1]Supported by an Australian Research Council Discovery grant.
[2]Supported by the National Natural Science Foundation of China and a Leverhulme Trust research grant.

[1], the literature has been dominated by conditional models, in particular, with the use of Markov fields and Markov chain Monte Carlo techniques. Another large branch of literature, mainly on irregularly positioned data, though, is concerned with the various methods of kriging which are based on parametric assumptions; see, for example, [6], Chapters 2–5.

In time series and regression, nonparametric methods have been very popular both for prediction and characterizing nonlinear dependence. No such development has taken place for spatial lattice models. Since the data are already on a grid, unless there are missing data, the prediction issue is less relevant, but there is still a need to explore and characterize nonlinear dependence relations. A rather obvious reason for the lack of progress is the curse of dimensionality. For a time series $\{Y_t\}$, a nonparametric regression $E[Y_t|Y_{t-1} = y]$ of $Y_t$ on its immediate predecessor is one-dimensional, and the corresponding Nadaraya–Watson (NW) estimator has good statistical properties. For spatial data $\{Y_{ij}\}$ on a grid, however, the conditional mean of $Y_{ij}$ given its closest neighbors $Y_{i-1,j}$, $Y_{i,j-1}$, $Y_{i+1,j}$ and $Y_{i,j+1}$ involves a four-dimensional nonparametric regression. Formally this can be carried out using the NW estimator, and an asymptotic theory can be constructed. In practice, however, this cannot be recommended unless the number of data points is extremely large.

In spite of these difficulties, there has been some recent theoretical work in this area. Kernel and nearest neighbor density estimates have been analyzed by Tran [33] and Tran and Yakowitz [34] under spatial mixing conditions. Clearly, in the marginal density estimation case, the curse of dimensionality is not an obstacle. The $L_1$ theory was established by Carbon, Hallin and Tran [4], and developed further by Hallin, Lu and Tran [15] under spatial stability conditions, including spatial linear and nonlinear processes, without imposing the less verifiable mixing conditions. The asymptotic normality of the kernel density estimator was also established for spatial linear processes by Hallin, Lu and Tran [14]. Finally, the NW kernel method and the local linear spatial conditional regressor were treated by Lu and Chen [21, 22], Hallin, Lu and Tran [16] and others. We have found these papers useful in developing our theory, but our perspective is rather different.

There are several ways of circumventing the curse of dimensionality in nonspatial regression. Perhaps the two most commonly used are semiparametric models, which in this context will be taken to mean partially linear models, and additive models. Actually, Cressie ([6], page 283) points out the possibility of trying such models for spatial data, noting that the nonlinear krige technique called disjunctive kriging (cf. [29]) takes as its starting point an additive decomposition. The problem, as seen from a traditional Markov field point of view, is that additivity clashes with the spatial Markov assumption. This is very different from the time series case where the partial linear autoregressive model (see [9])

$$Y_t = \beta Y_{t-1} + g(Y_{t-2}) + e_t$$

is a Markov model of second order if $\{e_t\}$ consists of independent and identically distributed (i.i.d.) random errors independent of $\{Y_{t-s}, s > 0\}$.

In the spatial case so far we have not been able to construct nonlinear additive or semiparametric models which are at the same time Markov. The problem can be illustrated by considering the line process $\{Y_i\}$. Assuming $\{Y_i\}$ to be Markov on the line and conditional Gaussian with density

$$p(y_i|y_{i-1}, y_{i+1}) = \frac{1}{\sqrt{2\pi}\sigma} e^{-(y_i - g(y_{i-1}) - h(y_{i+1}))^2/(2\sigma^2)},$$

it is easily seen using formulae (2.2) and (3.3) of [1] that the Markov field property implies $g(y) \equiv h(y) \equiv ay + b$ for two constants $a$ and $b$.

In ordinary regression, semiparametric and additive fitting can be thought of as an approximation of conditional quantities such as $E[Y_t|Y_{t-1}, \ldots, Y_{t-k}]$, and sometimes [31] interaction terms are included to improve this approximation. The approximation interpretation continues to be valid in the spatial case, so that semiparametric and additive models can be viewed as approximations to conditional expressions such as $E[Y_{ij}|Y_{i-1,j}, Y_{i,j-1}, Y_{i+1,j}, Y_{i,j+1}]$. The conditional spirit of Besag [1] is retained, being in terms of conditional means, however, rather than conditional probabilities. (Note that, also, in nonlinear time series, dependence is described by taking the conditional mean as a starting point; see, in particular, the contributions by Bjerve and Doksum [2] and Jones and Koch [19].) The conditional mean $E[Y_{ij}|Y_{i-1,j}, Y_{i,j-1}, Y_{i+1,j}, Y_{i,j+1}]$, say, is meaningful if first-order moments exist and if the conditional mean structure is invariant to spatial translations. Mathematically, the approximation consists in projecting this function on the set of semiparametric or additive functions. It is not claimed that there is a Markov field model, or any other conditional model, that can be exactly represented by this approximation. In this respect the situation is the same as for nonlinear disjunctive kriging, where the conditional mean of $Y_{ij}$ at a certain location is sought to be approximated by an additive decomposition going over all of the remaining observations (cf. [6], page 279). Classes of lattice models where there does exist an exact representation are the class of auto-Gaussian models (cf. [1]) or unilateral one-quadrant representations where $Y_{ij}$ is represented additively in terms of, say, $Y_{i-1,j}, Y_{i,j-1}$ only and an independent residual term (cf. [23]). But the former is linear, and the latter a "causal" unilateral expansion which may not be too realistic. In general, in the nonlinear spatial case, one must live with the approximative aspect. In practical time series modeling this is also the case, but in that situation at least one is able to write up a fairly general and exact model, where $Y$ can be expressed as an additive function of past values and an independent residual term. Fortunately, the asymptotic theory does not require the existence of such a representation.

The purpose of this paper is then to develop estimators for a spatial semiparametric (partially linear) structure and to derive their asymptotic properties. In the companion paper by Lu et al. [23], the additive approximation is analyzed using a different setup and different techniques of estimation. An advantage of using the partially linear approach is that a priori information concerning possible linearity of some of the components can be included in the model. More specifically, we will look at approximating the conditional mean function $m(X_{ij}, Z_{ij}) = E(Y_{ij} | X_{ij}, Z_{ij})$ by a semiparametric (partially linear) function of the form

$$(1.1) \qquad m_0(X_{ij}, Z_{ij}) = \mu + Z_{ij}^\tau \beta + g(X_{ij}),$$

such that $E[Y_{ij} - m_0(X_{ij}, Z_{ij})]^2$ or, equivalently, $E[m(X_{ij}, Z_{ij}) - m_0(X_{ij}, Z_{ij})]^2$ is minimized over a class of semiparametric functions of the form $m_0(X_{ij}, Z_{ij})$, subject to $E[g(X_{ij})] = 0$ for the identifiability of $m_0(X_{ij}, Z_{ij})$, where $\mu$ is an unknown parameter, $\beta = (\beta_1, \ldots, \beta_q)^\tau$ is a vector of unknown parameters, $g(\cdot)$ is an unknown function over $\mathbb{R}^p$, $Z_{ij} = (Z_{ij}^{(1)}, \ldots, Z_{ij}^{(q)})^\tau$ and $X_{ij} = (X_{ij}^{(1)}, \ldots, X_{ij}^{(p)})^\tau$ may contain both exogenous and endogenous variables, that is, neighboring values of $Y_{ij}$. Moreover, a component $Z_{ij}^{(r)}$ of $Z_{ij}$ or a component $X_{ij}^{(s)}$ of $X_{ij}$ may itself be a linear combination of neighboring values of $Y_{ij}$, as will be seen in Section 4, where $Z_{ij}^{(1)} = Y_{i-1,j} + Y_{i+1,j}$ and $X_{ij}^{(1)} = Y_{i,j-1} + Y_{i,j+1}$.

Motivation for using the form (1.1) for nonspatial data analysis can be found in [17]. As for the nonspatial case, estimating $g(\cdot)$ in model (1.1) may suffer from the curse of dimensionality when $g(\cdot)$ is not necessarily additive and $p \geq 3$. Thus, we will propose approximating $g(\cdot)$ by $g_a(\cdot)$, an additive marginal integration projector as detailed in Section 2 below. When $g(\cdot)$ itself is additive, that is, $g(x) = \sum_{l=1}^{p} g_l(x_l)$, $m_0(X_{ij}, Z_{ij})$ of (1.1) can be written as

$$(1.2) \qquad m_0(X_{ij}, Z_{ij}) = \mu + Z_{ij}^\tau \beta + \sum_{l=1}^{p} g_l(X_{ij}^{(l)}),$$

subject to $E[g_l(X_{ij}^{(l)})] = 0$ for all $1 \leq l \leq p$ for the identifiability of $m_0(X_{ij}, Z_{ij})$ in (1.2), where $g_l(\cdot)$, $l = 1, \ldots, p$, are all unknown one-dimensional functions over $\mathbb{R}^1$.

Our method of estimating $g(\cdot)$ or $g_a(\cdot)$ is based on an additive marginal integration projection on the set of additive functions, but where, unlike the backfitting case, the projection is taken with the product measure of $X_{ij}^{(l)}$ for $l = 1, \ldots, p$ (cf. [27]). This contrasts with the smoothed backfitting approach of Lu et al. [23], who base their work on an extension of the techniques of Mammen, Linton and Nelson [24] to the nonparametric spatial regression case. Marginal integration, although inferior to backfitting in asymptotic efficiency for purely additive models, seems well suited to the framework of partially linear estimation. In fact, in previous work (cf. [8]) in the independent regression case marginal integration has been used, and we do not know of any work extending the backfitting theory to the partially linear case. Marginal integration techniques are also applicable to the case where interactions are allowed between the $X_{ij}^{(k)}$-variables (cf. also the use of marginal integration for estimating interactions in ordinary regression problems).

We believe that our approach to analyzing spatial data is flexible. It permits nonlinearity and non-Gaussianity of real data. For example, reanalyzing the classical Mercer and Hall [26] wheat data set, one directional component appears to be nonlinear, and the fit is improved relative to earlier fits that have been linear. The presence of spatial dependence creates a host of new problems and, in particular, it has important effects on the estimation of the parametric component with asymptotic formulae different from those in the time series case.

The organization of the paper is as follows. Section 2 develops the kernel based marginal integration estimation procedure for the forms (1.1) and (1.2). Asymptotic properties of the proposed procedures are given in Section 3. Section 4 discusses an application of the proposed procedures to the Mercer and Hall data. A short conclusion is given in Section 5. Mathematical details are relegated to the Appendix.

## 2. Notation and definition of estimators.

As mentioned after (1.1), we are approximating the mean function $m(X_{ij}, Z_{ij}) = E[Y_{ij} | X_{ij}, Z_{ij}]$ by minimizing

$$E[Y_{ij} - m_0(X_{ij}, Z_{ij})]^2 = E[Y_{ij} - \mu - Z_{ij}^\tau \beta - g(X_{ij})]^2$$

over a class of semiparametric functions of the form $m_0(X_{ij}, Z_{ij}) = \mu + Z_{ij}^\tau \beta + g(X_{ij})$ with $E[g(X_{ij})] = 0$. Such a minimization problem is equivalent to minimizing

$$E[Y_{ij} - \mu - Z_{ij}^\tau \beta - g(X_{ij})]^2 = E[E\{(Y_{ij} - \mu - Z_{ij}^\tau \beta - g(X_{ij}))^2 | X_{ij}\}]$$

over some $(\mu, \beta, g)$. This implies that $g(X_{ij}) = E[(Y_{ij} - \mu - Z_{ij}^\tau \beta) | X_{ij}]$ and $\mu = E[Y_{ij} - Z_{ij}^\tau \beta]$, and $\beta$ is given by

$$\beta = (E[(Z_{ij} - E[Z_{ij} | X_{ij}])(Z_{ij} - E[Z_{ij} | X_{ij}])^\tau])^{-1}$$
$$\times E[(Z_{ij} - E[Z_{ij} | X_{ij}])(Y_{ij} - E[Y_{ij} | X_{ij}])],$$

provided that the inverse exists. This also shows that $m_0(X_{ij}, Z_{ij})$ is identifiable under the assumption of $E[g(X_{ij})] = 0$.

We now turn to estimation assuming that the data are available for $(Y_{ij}, X_{ij}, Z_{ij})$ for $1 \leq i \leq m, 1 \leq j \leq n$. Since nonparametric estimation is not much used for lattice data, and since the definitions of the estimators to be used later are quite involved notationally, we start by outlining the main steps in establishing estimators for $\mu$, $\beta$ and $g(\cdot)$ in (1.1) and then $g_l(\cdot), l = 1, 2, \ldots, p$, in (1.2). In the following, we give our outline in three steps.

*Step* 1. Estimating $\mu$ and $g(\cdot)$ assuming $\beta$ to be known.

For each fixed $\beta$, since $\mu = E[Y_{ij}] - E[Z_{ij}^{\tau} \beta] = \mu_Y - \mu_Z^{\tau} \beta$, $\mu$ can be estimated by $\hat{\mu}(\beta) = \overline{Y} - \overline{Z}^{\tau} \beta$, where $\mu_Y = E[Y_{ij}]$, $\mu_Z = (\mu_Z^{(1)}, \ldots, \mu_Z^{(q)})^{\tau} = E[Z_{ij}]$, $\overline{Y} = \frac{1}{mn} \sum_{i=1}^{m} \sum_{j=1}^{n} Y_{ij}$ and $\overline{Z} = \frac{1}{mn} \sum_{i=1}^{m} \sum_{j=1}^{n} Z_{ij}$.

Moreover, the conditional expectation

$$
\begin{aligned}
g(x) = g(x, \beta) &= E[(Y_{ij} - \mu - Z_{ij}^{\tau} \beta) | X_{ij} = x] \\
&= E[(Y_{ij} - E[Y_{ij}] - (Z_{ij} - E[Z_{ij}])^{\tau} \beta) | X_{ij} = x]
\end{aligned}
$$

can be estimated by standard local linear estimation ([7], page 19), with $\hat{g}_{m,n}(x, \beta) = \hat{a}_0(\beta)$ satisfying

$$
\begin{aligned}
(2.1) \quad & (\hat{a}_0(\beta), \hat{a}_1(\beta)) \\
&= \underset{(a_0, a_1) \in \mathbb{R}^1 \times \mathbb{R}^p}{\arg \min} \sum_{i=1}^{m} \sum_{j=1}^{n} (\tilde{Y}_{ij} - \tilde{Z}_{ij}^{\tau} \beta - a_0 - a_1^{\tau}(X_{ij} - x))^2 K_{ij}(x, b),
\end{aligned}
$$

where $\tilde{Y}_{ij} = Y_{ij} - \overline{Y}$ and $\tilde{Z}_{ij} = (\tilde{Z}_{ij}^{(1)}, \ldots, \tilde{Z}_{ij}^{(q)})^{\tau} = Z_{ij} - \overline{Z}$.

*Step* 2. Marginal integration to obtain $g_1, \ldots, g_p$ of (1.2).

The idea of the marginal integration estimator is best explained if $g(\cdot)$ is itself additive, that is, if

$$
g(X_{ij}) = g(X_{ij}^{(1)}, \ldots, X_{ij}^{(p)}) = \sum_{l=1}^{p} g_l(X_{ij}^{(l)}).
$$

Then, since $E[g_l(X_{ij}^{(l)})] = 0$ for $l = 1, \ldots, p$, for $k$ fixed,

$$
g_k(x_k) = E[g(X_{ij}^{(1)}, \ldots, x_k, \ldots, X_{ij}^{(p)})]
$$

and an estimate of $g_k$ is obtained by keeping $X_{ij}^{(k)}$ fixed at $x_k$ and then taking the average over the remaining variables $X_{ij}^{(1)}, \ldots, X_{ij}^{(k-1)}, X_{ij}^{(k+1)}, \ldots, X_{ij}^{(p)}$. This marginal integration operation can be implemented irrespective of whether or not $g(\cdot)$ is additive. If the additivity does not hold, as mentioned in the Introduction, the marginal integration amounts to a projection on the space of additive functions of $X_{ij}^{(l)}$, $l = 1, \ldots, p$, taken with respect to the product measure of $X_{ij}^{(l)}, l = 1, \ldots, p$, obtaining the approximation $g_a(x, \beta) = \sum_{l=1}^{p} P_{l,\omega}(X_{ij}^{(l)}, \beta)$, which will be detailed below with $\beta$ appearing linearly in the expression. In addition, it has been found convenient to introduce a pair of weight functions $(w_k, w_{(-k)})$ in the estimation of each component, hence, the index $w$ in $P_{l,w}$. The details are given in (2.7)–(2.9) below.

*Step* 3. Estimating $\beta$.

The last step consists in estimating $\beta$. This is done by weighted least squares, and it is easy since $\beta$ enters linearly in our expressions. In fact, using the expression of $g(x, \beta)$ in step 1, we obtain the weighted least squares estimator $\hat{\beta}$ of $\beta$ in (2.10) below. Finally, this is re-introduced in the expressions for $\hat{\mu}$ and $\hat{P}$ resulting in the estimates in (2.11) and (2.12) below. In the following, steps 1–3 are written correspondingly in more detail.

*Step* 1. To write our expression for $(\hat{a}_0(\beta), \hat{a}_1(\beta))$ in (2.1), we need to introduce some more notation. Let $K_{ij} = K_{ij}(x, b) = \prod_{l=1}^{p} K(\frac{X_{ij}^{(l)} - x_l}{b_l})$, with $b = b_{m,n} = (b_1, \ldots, b_p)$, $b_l = b_{l,m,n}$ being a sequence of bandwidths for the $l$th covariate variable $X_{ij}^{(l)}$, tending to zero as $(m, n)$ tends to infinity, and $K(\cdot)$ is a bounded kernel function on $\mathbb{R}^1$ (when we do the asymptotic analysis in Section 3, we need to introduce a more refined choice of bandwidths, as is explained just before stating Assumption 3.6). Denote

$$\mathcal{X}_{ij} = \mathcal{X}_{ij}(x, b) = \left( \frac{(X_{ij}^{(1)} - x_1)}{b_1}, \ldots, \frac{(X_{ij}^{(p)} - x_p)}{b_p} \right)^\tau,$$

and let $b_\pi = \prod_{l=1}^{p} b_l$. We define

$$u_{m,n,l_1 l_2} = (mnb_\pi)^{-1} \sum_{i=1}^{m} \sum_{j=1}^{n} (\mathcal{X}_{ij}(x, b))_{l_1} (\mathcal{X}_{ij}(x, b))_{l_2} K_{ij}(x, b),$$

$$0 \leq l_1, l_2 \leq p,$$

where $(\mathcal{X}_{ij}(x, b))_l = (X_{ij}^{(l)} - x_l)/b_l$ for $1 \leq l \leq p$. We then let $(\mathcal{X}_{ij}(x, b))_0 \equiv 1$ and define

$$v_{m,n,l}(\beta) = (mnb_\pi)^{-1} \sum_{i=1}^{m} \sum_{j=1}^{n} (\tilde{Y}_{ij} - \tilde{Z}_{ij}^\tau \beta)(\mathcal{X}_{ij}(x, b))_l K_{ij}(x, b) \tag{2.3}$$

and where, as before, $\tilde{Y}_{ij} = Y_{ij} - \bar{Y}$ and $\tilde{Z}_{ij} = Z_{ij} - \bar{Z}$.

Note that $v_{m,n,l}(\beta)$ can be decomposed as

$$(2.4) \qquad v_{m,n,l}(\beta) = v_{m,n,l}^{(0)} - \sum_{s=1}^{q} \beta_s v_{m,n,l}^{(s)} \qquad \text{for } l = 0, 1, \ldots, p,$$

in which

$$v_{m,n,l}^{(0)} = v_{m,n,l}^{(0)}(x, b)$$
$$= (mnb_\pi)^{-1} \sum_{i=1}^{m} \sum_{j=1}^{n} \tilde{Y}_{ij} \left( \mathcal{X}_{ij}(x, b) \right)_l K_{ij}(x, b),$$

$$v_{m,n,l}^{(s)} = v_{m,n,l}^{(s)}(x, b)$$
$$= (mnb_\pi)^{-1} \sum_{i=1}^{m} \sum_{j=1}^{n} \tilde{Z}_{ij}^{(s)} (\mathcal{X}_{ij}(x, b))_l K_{ij}(x, b), \qquad 1 \le s \le q.$$

We can then express the local linear estimates in (2.1) as

$$(2.5) \qquad (\hat{a}_0(\beta), \hat{a}_1(\beta) \odot b)^\tau = U_{m,n}^{-1} V_{m,n}(\beta),$$

where $\odot$ is the operation of the component-wise product, that is, $a_1 \odot b = (a_{11}b_1, \ldots, a_{1p}b_p)$ for $a_1 = (a_{11}, \ldots, a_{1p})$ and $b = (b_1, \ldots, b_p)$,

$$(2.6) \qquad V_{m,n}(\beta) = \begin{pmatrix} v_{m,n,0}(\beta) \\ V_{m,n,1}(\beta) \end{pmatrix}, \qquad U_{m,n} = \begin{pmatrix} u_{m,n,00} & U_{m,n,01} \\ U_{m,n,10} & U_{m,n,11} \end{pmatrix},$$

where $U_{m,n,10} = U_{m,n,01}^\tau = (u_{m,n,01}, \ldots, u_{m,n,0p})^\tau$ and $U_{m,n,11}$ is the $p \times p$ matrix defined by $u_{m,n,l_1 l_2}$, with $l_1, l_2 = 1, \ldots, p$, in (2.2). Moreover, $V_{m,n,1}(\beta) = (v_{m,n,1}(\beta), \ldots, v_{m,n,p}(\beta))^\tau$, with $v_{m,n,l}(\beta)$ as defined in (2.3). Analogously for $V_{m,n}$, we may define $V_{m,n}^{(0)}$ and $V_{m,n}^{(s)}$ in terms of $v_{m,n}^{(0)}$ and $v_{m,n}^{(s)}$. Then taking the first component with $\gamma = (1, 0, \ldots, 0)^\tau \in \mathbb{R}^{1+p}$,

$$\hat{g}_{m,n}(x, \beta) = \gamma^\tau U_{m,n}^{-1}(x) V_{m,n}(x, \beta)$$
$$= \gamma^\tau U_{m,n}^{-1}(x) V_{m,n}^{(0)}(x) - \sum_{s=1}^{q} \beta_s \gamma^\tau U_{m,n}^{-1}(x) V_{m,n}^{(s)}(x)$$
$$= H_{m,n}^{(0)}(x) - \beta^\tau H_{m,n}(x),$$

where $H_{m,n}(x) = (H_{m,n}^{(1)}(x), \ldots, H_{m,n}^{(q)}(x))^\tau$, with $H_{m,n}^{(s)}(x) = \gamma^\tau U_{m,n}^{-1}(x) V_{m,n}^{(s)}(x)$, $1 \le s \le q$. Clearly, $H_{m,n}^{(s)}(x)$ is the local linear estimator of $H^{(s)}(x) = E[(Z_{ij}^{(s)} - \mu_Z^{(s)})|X_{ij} = x]$, $1 \le s \le q$.

We now define $Z_{ij}^{(0)} = Y_{ij}$ and $\mu_Z^{(0)} = \mu_Y$ such that $H^{(0)}(x) = E[(Z_{ij}^{(0)} - \mu_Z^{(0)})|X_{ij} = x] = E[Y_{ij} - \mu_Y|X_{ij} = x]$ and $H(x) = (H^{(1)}(x), \ldots, H^{(q)}(x))^\tau = E[(Z_{ij} - \mu_Z)|X_{ij} = x]$. It follows that $g(x, \beta) = H^{(0)}(x) - \beta^\tau H(x)$, which equals $g(x)$ under (1.1) irrespective of whether $g$ itself is additive.

*Step* 2. Let $w_{(-k)}(\cdot)$ be a weight function defined on $\mathbb{R}^{p-1}$ such that $E[w_{(-k)}(X_{ij}^{(-k)})] = 1$, and $w_k(x_k) = I_{[-L_k, L_k]}(x_k)$ defined on $\mathbb{R}^1$ for some large $L_k > 0$, with

$$X_{ij}^{(-k)} = (X_{ij}^{(1)}, \ldots, X_{ij}^{(k-1)}, X_{ij}^{(k+1)}, \ldots, X_{ij}^{(p)}),$$

where $I_A(x)$ is the conventional indicator function.

For a given $\beta$, consider the marginal projection

$$
\begin{align}
(2.7) \quad P_{k,w}(x_k, \beta) = E[g(X_{ij}^{(1)}, \ldots, X_{ij}^{(k-1)}, x_k, \\
X_{ij}^{(k+1)}, \ldots, X_{ij}^{(p)}, \beta) w_{(-k)}(X_{ij}^{(-k)})] w_k(x_k).
\end{align}
$$

It is easily seen that if $g$ is additive as in (1.2), then, for $-L_k \leq x_k \leq L_k$, $P_{k,w}(x_k, \beta) = g_k(x_k)$ up to a constant since it is assumed that $E[w_{(-k)}(X_{ij}^{(-k)})] = 1$. In general, $g_a(x, \beta) = \sum_{l=1}^{p} P_{l,w}(x_l, \beta)$ is an additive marginal projection approximation to $g(x)$ in (1.1) up to a constant in the region $x \in \prod_{l=1}^{p}[-L_l, L_l]$. The quantity $P_{k,w}(x_k, \beta)$ can then be estimated by the spatial locally linear marginal integration estimator

$$
\begin{align}
\widehat{P}_{k,w}(x_k, \beta) = (mn)^{-1} \sum_{i=1}^{m} \sum_{j=1}^{n} \hat{g}_{m,n}(X_{ij}^{(1)}, \ldots, X_{ij}^{(k-1)}, x_k, \\
(2.8) \quad\quad X_{ij}^{(k+1)}, \ldots, X_{ij}^{(p)}, \beta) w_{(-k)}(X_{ij}^{(-k)}) w_k(x_k) \\
= \hat{P}_{k,w}^{(0)}(x_k) - \sum_{s=1}^{q} \beta_s \hat{P}_{k,w}^{(s)}(x_k) = \hat{P}_{k,w}^{(0)}(x_k) - \beta^\tau \hat{P}_{k,w}^Z(x_k),
\end{align}
$$

where

$$
\begin{align}
\hat{P}_{k,w}^{(s)}(x_k) = \frac{1}{mn} \sum_{i=1}^{m} \sum_{j=1}^{n} H_{m,n}^{(s)}(X_{ij}^{(1)}, \ldots, X_{ij}^{(k-1)}, x_k, \\
X_{ij}^{(k+1)}, \ldots, X_{ij}^{(p)}) w_{(-k)}(X_{ij}^{(-k)}) w_k(x_k)
\end{align}
$$

is the estimator of

$$
\begin{align}
P_{k,w}^{(s)}(x_k) = E[H^{(s)}(X_{ij}^{(1)}, \ldots, X_{ij}^{(k-1)}, x_k, \\
X_{ij}^{(k+1)}, \ldots, X_{ij}^{(p)}) w_{(-k)}(X_{ij}^{(-k)})] w_k(x_k),
\end{align}
$$

for $0 \leq s \leq q$, and $P_{k,w}^Z(x_k) = (P_{k,w}^{(1)}(x_k), \ldots, P_{k,w}^{(q)}(x_k))^\tau$ is estimated by

$$\hat{P}_{k,w}^Z(x_k) = (\hat{P}_{k,w}^{(1)}(x_k), \ldots, \hat{P}_{k,w}^{(q)}(x_k))^\tau.$$

Here, we add the weight function $w_k(x_k) = I_{[-L_k, L_k]}(x_k)$ in the definition of $\hat{P}_{k,w}^{(s)}(x_k)$, since we are only interested in the points of $x_k \in [-L_k, L_k]$ for some large $L_k$. In practice, we may use a sample centered version of $\hat{P}_{k,w}^{(s)}(x_k)$ as the estimator of $P_{k,w}^{(s)}(x_k)$. Clearly, we have $P_{k,w}(x_k, \beta) = P_{k,w}^{(0)}(x_k) - \beta^\tau P_{k,w}^Z(x_k)$. Thus, for every $\beta$, $g(x) = g(x, \beta)$ of (1.1) [or rather the approximation $g_a(x, \beta)$ if (1.2) does not hold] can be estimated by

$$(2.9) \qquad \widehat{\widehat{g}}(x, \beta) = \sum_{l=1}^p \hat{P}_{l,w}(x_l, \beta) = \sum_{l=1}^p \hat{P}_{l,w}^{(0)}(x_l) - \beta^\tau \sum_{l=1}^p \hat{P}_{l,w}^Z(x_l).$$

*Step* 3.  We can finally obtain the least squares estimator of $\beta$ by

$$\hat{\beta} = \arg\min_{\beta \in \mathbb{R}^q} \sum_{i=1}^m \sum_{j=1}^n (\tilde{Y}_{ij} - \tilde{Z}_{ij}^\tau \beta - \widehat{\widehat{g}}(X_{ij}, \beta))^2$$

$$(2.10)$$

$$= \arg\min_{\beta \in \mathbb{R}^q} \sum_{i=1}^m \sum_{j=1}^n (\hat{Y}_{ij}^* - (\widehat{Z}_{ij}^*)^\tau \beta)^2,$$

where $\hat{Y}_{ij}^* = \tilde{Y}_{ij} - \sum_{l=1}^p \hat{P}_{l,w}^{(0)}(X_{ij}^{(l)})$ and $\widehat{Z}_{ij}^* = \tilde{Z}_{ij} - \sum_{l=1}^p \hat{P}_{l,w}^Z(X_{ij}^{(l)})$. Therefore,

$$(2.11) \quad \hat{\beta} = \left( \sum_{i=1}^m \sum_{j=1}^n \widehat{Z}_{ij}^* (\widehat{Z}_{ij}^*)^\tau \right)^{-1} \left( \sum_{i=1}^m \sum_{j=1}^n \hat{Y}_{ij}^* \widehat{Z}_{ij}^* \right) \quad \text{and} \quad \hat{\mu} = \overline{Y} - \hat{\beta}^\tau \overline{Z}.$$

We then insert $\hat{\beta}$ in $\hat{a}_0(\beta) = \hat{g}_{m,n}(x, \beta)$ to obtain $\hat{a}_0(\hat{\beta}) = \hat{g}_{m,n}(x, \hat{\beta})$. In view of this, the spatial local linear projection estimator of $P_k(x_k)$ can be defined by

$$\widehat{\widehat{P}}_{k,w}(x_k) = (mn)^{-1} \sum_{i=1}^m \sum_{j=1}^n \hat{g}_{m,n}(X_{ij}^{(1)}, \ldots, X_{ij}^{(k-1)}, x_k,$$

$$(2.12)$$

$$X_{ij}^{(k+1)}, \ldots, X_{ij}^{(p)}; \hat{\beta}) w_{(-k)}(X_{ij}^{(-k)}),$$

and for $x_k \in [-L_k, L_k]$, this would estimate $g_k(x_k)$ up to a constant when (1.2) holds. To ensure $E[g_k(X_{ij}^{(k)})] = 0$, we may rewrite $\widehat{\widehat{P}}_{k,w}(x_k) - \hat{\mu}_P(k)$ for the estimate of $g_k(x_k)$ in (1.2), where $\hat{\mu}_P(k) = \frac{1}{mn} \sum_{i=1}^m \sum_{j=1}^n \widehat{\widehat{P}}_{k,w}(X_{ij}^{(k)})$.

For the least squares estimator, $\hat{\beta}$, and $\widehat{\widehat{P}}_{k,w}(\cdot)$, we establish some asymptotic distributions under mild conditions in Section 3.

**3. Asymptotic properties.** Let $\mathcal{I}_{m,n}$ be the rectangular region defined by $\mathcal{I}_{m,n} = \{(i,j) : i, j \in \mathbb{Z}^2, 1 \le i \le m, 1 \le j \le n\}$. We observe $\{(Y_{ij}, X_{ij}, Z_{ij})\}$ on $\mathcal{I}_{m,n}$ with a sample size of $mn$.

In this paper we write $(m, n) \to \infty$ if

$$(3.1) \qquad\qquad\qquad \min\{m, n\} \to \infty.$$

In [33] it is required, in addition, that $m$ and $n$ tend to infinity at the same rate:

$$(3.2) \qquad C_1 < |m/n| < C_2 \qquad \text{for some } 0 < C_1 < C_2 < \infty.$$

Let $\{(Y_{ij}, X_{ij}, Z_{ij})\}$ be a strictly stationary random field indexed by $(i,j) \in \mathbb{Z}^2$. A point $(i,j)$ in $\mathbb{Z}^2$ is referred to as a site. Let $S$ and $S'$ be two sets of sites. The Borel fields $\mathcal{B}(S) = \mathcal{B}(Y_{ij}, X_{ij}, Z_{ij}, (i,j) \in S)$ and $\mathcal{B}(S') = \mathcal{B}(Y_{ij}, X_{ij}, Z_{ij}, (i,j) \in S')$ are the $\sigma$-fields generated by the random variables $(Y_{ij}, X_{ij}, Z_{ij})$, with $(i,j)$ being elements of $S$ and $S'$, respectively. We will assume that the variables $(Y_{ij}, X_{ij}, Z_{ij})$ satisfy the following mixing condition (cf. [33]): There exists a function $\varphi(t) \downarrow 0$ as $t \to \infty$, such that, whenever $S, S' \subset \mathbb{Z}^2$,

$$
\begin{aligned}
(3.3) \quad \alpha(\mathcal{B}(S), \mathcal{B}(S')) &= \sup_{\{A \in \mathcal{B}(S), B \in \mathcal{B}(S')\}} \{|P(AB) - P(A)P(B)|\} \\
&\leq \tilde{f}(\mathrm{Card}(S), \mathrm{Card}(S')) \varphi(\tilde{d}(S, S')),
\end{aligned}
$$

where $\mathrm{Card}(S)$ denotes the cardinality of $S$, and $\tilde{d}$ is the distance defined by

$$\tilde{d}(S, S') = \min\{\sqrt{|i - i'|^2 + |j - j'|^2} : (i,j) \in S, (i',j') \in S'\}.$$

Here $\tilde{f}$ is a symmetric positive function nondecreasing in each variable. Throughout the paper, we only assume that $\tilde{f}$ satisfies

$$(3.4) \qquad \tilde{f}(n, m) \leq \min\{m, n\}.$$

If $\tilde{f} \equiv 1$, then the spatial process $\{(Y_{ij}, X_{ij}, Z_{ij})\}$ is called strongly mixing. Condition (3.4) holds in many cases. Examples can be found in [30]. For relevant work on random fields, see, for example, [3, 5, 12, 13, 20, 28, 32, 35].

To state and prove our main results, we introduce the following assumptions.

ASSUMPTION 3.1. Assume that the process $\{(Y_{ij}, X_{ij}, Z_{ij}) : (i,j) \in \mathbb{Z}^2\}$ is strictly stationary. The joint probability density $f_s(x_1, \ldots, x_s)$ of $(X_{i_1 j_1}, \ldots, X_{i_s j_s})$ exists and is bounded for $s = 1, \ldots, 2r-1$, where $r$ is some positive integer such that Assumption 3.2(ii) below holds. For $s = 1$, we write $f(x)$ for $f_1(x_1)$, the density function of $X_{ij}$.

ASSUMPTION 3.2. (i) Let $Z_{ij}^* = Z_{ij} - \mu_Z - \sum_{l=1}^p P_{l,w}^Z(X_{ij}^{(l)})$ and $B^{ZZ} = E[Z_{11}^*(Z_{11}^*)^\tau]$. The inverse matrix of $B^{ZZ}$ exists. Let $Y_{ij}^* = Y_{ij} - \mu_Y - \sum_{l=1}^p P_{l,w}^{(0)} \times (X_{ij}^{(l)})$ and $R_{ij} = Z_{ij}^*(Y_{ij}^* - Z_{ij}^{*\tau}\beta)$. Assume that the matrix $\Sigma_B = \sum_{i=-\infty}^{\infty} \times \sum_{j=-\infty}^{\infty} E[(R_{00} - \mu_B)(R_{ij} - \mu_B)^\tau]$ is finite.

(ii) Suppose there is some $\lambda > 2$ such that $E[|Y_{ij}|^{\lambda r}] < \infty$ for $r$ as defined in Assumption 3.1.

Assumption 3.3.   The mixing coefficient $\varphi$ defined in (3.3) satisfies

$$\lim_{T \to \infty} T^a \sum_{t=T}^{\infty} t^{2r-1} \varphi(t)^{(\lambda r - 2)/(\lambda r)} = 0 \tag{3.5}$$

for some constant $a > \max(\frac{2(r\lambda+2)}{\lambda r}, \frac{2r(\lambda r - 2)}{2+\lambda r - 4r})$, with $\lambda > 4 - \frac{2}{r}$ as in Assumption 3.2(ii). In addition, the coefficient function $\tilde{f}$ involved in (3.3) satisfies (3.4).

Assumption 3.4.   (i) The functions $g(\cdot)$ in (1.1) and $g_l(\cdot)$ for $1 \le l \le p$ in (1.2) have bounded and continuous derivatives up to order 2. In addition, the function $g(\cdot)$ has a second-order derivative matrix $g''(\cdot)$ (of dimension $p \times p$), which is uniformly continuous on $\mathbb{R}^p$.

(ii) For each $k$, $1 \le k \le p$, the weight function $\{w_{(-k)}(\cdot)\}$ is uniformly continuous on $\mathbb{R}^{p-1}$ and bounded on the compact support $S_w^{(-k)}$ of $w_{(-k)}(\cdot)$. In addition, $E[w_{(-k)}(X_{ij}^{(-k)})] = 1$. Let $S_W = S_{W,k} = S_w^{(-k)} \times [-L_k, L_k]$ be the compact support of $W(x) = W(x^{(-k)}, x_k) = w_{(-k)}(x^{(-k)}) \cdot I_{[-L_k, L_k]}(x_k)$. In addition, let $\inf_{x \in S_W} f(x) > 0$ hold.

Assumption 3.5.   The function $K(x)$ is a symmetric and bounded probability density function on $\mathbb{R}^1$ with compact support, $C_K$, and finite variance such that $|K(x) - K(y)| \le M|x - y|$ for $x, y \in C_K$ and $0 < M < \infty$.

When we are estimating the marginal projector $P_k$, the bandwidth $b_k$ associated with this component has to tend to zero at a rate slower than $b_l$ for $l \ne k$. This means that, for each $k$, $1 \le k \le p$, we need a separate set of bandwidths $b_1^{(k)}, \ldots, b_p^{(k)}$ such that $b_k^{(k)}$ tends to zero slower than $b_l^{(k)}$ for all $l \ne k$. Correspondingly, we get $p$ different products $b_\pi^{(k)} = \prod_{l=1}^p b_l^{(k)}$. Since in the following we will analyze one component $\hat{P}_k$ at a time, to simplify notation we omit the superscript $(k)$ and write $b_k$, $b_l, l \ne k$, and $b_\pi$ instead of $b_k^{(k)}$, $b_l^{(k)}, l \ne k$, and $b_\pi^{(k)}$. It will be seen that this slight abuse of notation does not lead to interpretational difficulties in the proofs. To have consistency in notation, Assumptions 3.6 and 3.6′ below are also formulated using this notational simplification. Throughout the whole paper, we use $l$ as any arbitrary index, while leaving $k$ for the fixed and specified index as suggested by a referee.

Assumption 3.6.   (i) Let $b_\pi$ be as defined before. The bandwidths satisfy

$$\lim_{(m,n) \to \infty} \max_{1 \le l \le p} b_l = 0,$$

$$\lim_{(m,n)\to\infty} mnb_\pi^{1+2/r} = \infty,$$

$$\liminf_{(m,n)\to\infty} mnb_\pi^{2(r-1)a+2(\lambda r-2)/((a+2)\lambda)} > 0$$

for some integer $r \geq 3$ and some $\lambda > 2$ being the same as in Assumptions 3.1 and 3.2.

(ii) In addition, for some integer $r \geq 3$, the $k$th component satisfies

$$\limsup_{(m,n)\to\infty} mnb_k^5 < \infty,$$

$$\lim_{(m,n)\to\infty} \frac{\max_{1\leq l\neq k\leq p} b_l}{b_k} = 0,$$

$$\lim_{(m,n)\to\infty} mnb_k^{4(2+r)/(2r-1)} = \infty.$$

REMARK 3.1. (i) Assumptions 3.1, 3.2, 3.4 and 3.5 are relatively mild in this kind of problem, and can be justified in detail. For example, Assumption 3.1 is quite natural and corresponds to that used for the nonspatial case. Assumption 3.2(i) is necessary for the establishment of asymptotic normality in the semiparametric setting. As can be seen from Theorem 3.1 below, the condition on the existence of the inverse matrix, $(B^{ZZ})^{-1}$, is required in the formulation of that theorem. Moreover, Assumption 3.2(i) corresponds to those used for the nonspatial case. Assumption 3.2(ii) is needed as the existence of moments of higher than second order is required for this kind of problem when uniform convergence for nonparametric regression estimation is involved. Assumption 3.4(ii) is required due to the use of such a weight function. The continuity condition on the kernel function is quite natural and easily satisfied.

(ii) As for the nonspatial case (see Condition A of [8]), some technical conditions are needed when marginal integration techniques are employed. In addition, some other technical conditions are required for the spatial case. Condition (3.5) requires some kind of rate of convergence for the mixing coefficient. It holds automatically when the mixing coefficient decreases to zero exponentially. For the nonspatial case, similar conditions have been used. See, for example, Condition A(vi) of [8]. For the spatial case, Assumption 3.6 requires that, when one of the bandwidths is proportional to $(mn)^{-1/5}$, the optimal choice under a conventional criterion, the other bandwidths need to converge to zero with a rate related to $(mn)^{-1/5}$. Assumption 3.6 is quite complex in general. However, it holds in some cases. For example, when we choose $p = 2$, $r = 3$, $\lambda = 4$, $a = 31$, $k = 1$, $b_1 = (mn)^{-1/5}$ and $b_2 = (mn)^{-2/5+\eta}$ for some $0 < \eta < \frac{1}{5}$, both (i) and (ii) hold. For instance,

$$\liminf_{(m,n)\to\infty} mnb_\pi^{2(r-1)a+2(\lambda r-2)/((a+2)\lambda)}$$

$$= \liminf_{(m,n)\to\infty} (mn)^{(19/55)+(12/11)\eta} = \infty > 0$$

and

$$\lim_{(m,n)\to\infty} mnb_\pi^{1+2/r} = \lim_{(m,n)\to\infty} (mn)^{(5/3)\eta} = \infty.$$

(iii) Similarly to the nonspatial case ([8], Remark 10), we assume that all the nonparametric components are only two times continuously differentiable and, thus, the optimal bandwidth $b_k$ is proportional to $(mn)^{-1/5}$. As a result, Assumption 3.6 basically implies $p \le 4$. For our case, the assumption of $p \le 4$ is just sufficient for us to use an additive model to approximate the conditional mean $E[Y_{ij}|Y_{i-1,j}, Y_{i,j-1}, Y_{i+1,j}, Y_{i,j+1}]$ by $g_1(Y_{i-1,j}) + g_2(Y_{i,j-1}) + g_3(Y_{i+1,j}) + g_4(Y_{i,j+1})$, with each $g_i(\cdot)$ being an unknown function. In addition, for our case study in Section 4, we need only to use an additive model of the form $g_1(X_{ij}^{(1)}) + g_2(X_{ij}^{(2)})$ to approximate the conditional mean, where $X_{ij}^{(1)} = Y_{i,j-1} + Y_{i,j+1}$ and $X_{ij}^{(2)} = Y_{i-1,j} + Y_{i+1,j}$. Nevertheless, we may ensure that the marginal integration method still works for the case of $p \ge 5$ and achieves the optimal rate of convergence by using a high-order kernel of the form

(3.6)
$$\int K(x)\,dx = 1,$$
$$\int x^i K(x)\,dx = 0 \qquad \text{for } i = 1, \ldots, I-1 \quad \text{and}$$
$$\int x^I K(x) \ne 0$$

for $I \ge 2$, as discussed in [18] for the nonspatial case, where $I$ is the order of smoothness of the nonparametric components. To ensure that the conclusions of the main results hold for this case, we need to replace Assumptions 3.4–3.6 by Assumptions 3.4′–3.6′ below:

ASSUMPTION 3.4′.    (i) The functions $g(\cdot)$ in (1.1) and $g_l(\cdot)$ for $1 \le l \le p$ in (1.2) have bounded and continuous derivatives up to order $I \ge 2$. In addition, the function $g(\cdot)$ has an $I$-order derivative matrix $g^{(I)}(\cdot)$ (of dimension $p \times p \times \cdots \times p$) which is uniformly continuous on $\mathbb{R}^p$.

(ii) Assumption 3.4(ii) holds.

ASSUMPTION 3.5′.    Assumption 3.5(i) holds. In addition, the kernel function satisfies (3.6).

ASSUMPTION 3.6′.    (i) Assumption 3.6(i) holds.

(ii) In addition, for the $k$th component,

$$\limsup_{(m,n)\to\infty} mn b_k^{2I+1} < \infty,$$

$$\lim_{(m,n)\to\infty} \frac{\max_{1\le l\neq k\le p} b_l}{b_k} = 0,$$

$$\lim_{(m,n)\to\infty} mn b_k^{4(2+r)/(2r-1)} = \infty$$

for $\lambda > 2$ and some integer $r \ge 3$.

After Assumptions 3.4–3.6 are replaced by Assumptions 3.4′–3.6′, we may show that the conclusions of the results remain true. Under Assumptions 3.4′–3.6′, we will need to make changes at several places in the proofs of Lemmas A.3–A.5 and Theorems 3.1 and 3.2. Apart from replacing Assumptions 3.4–3.6 by Assumptions 3.4′–3.6′ in their conditions, we need to replace $\sum_{k=1}^{p} b_k^2$ by $\sum_{k=1}^{p} b_k^I$ and $\mu_2(K) = \int u^2 K(u)\,du$ by $\mu_I(K) = \int u^I K(u)\,du$, for example, in several relevant places.

To verify Assumption 3.6′, we can choose (remember the notational simplification introduced just before Assumption 3.6) the optimal bandwidth $b_k \sim (mn)^{-1/(2I+1)}$ and $b_l \sim (mn)^{-2/(2I+1)+\eta}$, with $0 < \eta < \frac{1}{2I+1}$ for all $l \neq k$. In this case, it is not difficult to verify Assumption 3.6′ for the case $p \ge 5$. As expected, the order of the smoothness $I$ needs to be greater than 2. For example, it is easy to see that Assumption 3.6′ holds for the case $p = 6$ when we choose $a = 31$, $r = 3$, $\lambda = 4$ and $I > 4 + \frac{1}{2}$. For instance, on the one hand, in order to make sure that the condition $\lim_{(m,n)\to\infty} \frac{\max_{1\le l\neq k\le p} b_l}{b_k} = 0$ holds, we need to have $0 < \eta < \frac{1}{2I+1}$. On the other hand, in order to ensure that

$$\liminf_{(m,n)\to\infty} mn b_\pi^{2(r-1)a+2(\lambda r-2)/((a+2)\lambda)}$$

$$= \liminf_{(m,n)\to\infty} (mn)^{(2I-11)/(2I+1)+(60/11)\eta} = \infty > 0$$

and

$$\lim_{(m,n)\to\infty} mn b_\pi^{1+2/r} = \lim_{(m,n)\to\infty} (mn)^{(6I-52)/(3(2I+1))+(25/3)\eta} = \infty$$

both hold, we need to assume $\eta > \frac{52-6I}{25(2I+1)}$. Thus, we can choose $\eta$ such that $\frac{52-6I}{25(2I+1)} < \eta < \frac{1}{2I+1}$ when $I > 4 + \frac{1}{2}$. The last equation of Assumption 3.6′(ii) holds automatically when $I > 4 + \frac{1}{2}$.

As pointed out by a referee, in general, to ensure that Assumption 3.6′ holds, we will need to choose $\eta$ such that $\frac{[2(p-1)+1](1+2/r)-(2I+1)}{(p-1)(1+2/r)} < \eta < \frac{1}{2I+1}$, which implies that $(I, p, r)$ does need to satisfy $I > \frac{(p-1)r+2p}{2r}$.

This suggests that, in order to achieve the rate-optimal property, we will need to allow that smoothness increases with dimensions. This is well known and has been used in some recent papers for the nonspatial case (see Conditions A5, A7 and NW2–NW3 of [18]).

(iv) Assumptions 3.2(ii), 3.3 and 3.6 together require the existence of $E[|Y_{ij}|^{10+\epsilon}]$ for some small $\epsilon > 0$. This may look like a strong moment condition. However, this is weaker than $E[|Y_{ij}|^k] < \infty$ for $k = 1, 2, \ldots$ and $E[e^{|Y_{ij}|}] < \infty$ corresponding to those used in the nonspatial case.

We can now state the asymptotic properties of the marginal integration estimators for both the parametric and nonparametric components. Recall that $Z_{ij}^* = Z_{ij} - \mu_Z - \sum_{l=1}^p P_{l,w}(X_{ij}^{(l)})$, $Y_{ij}^* = Y_{ij} - \mu_Y - \sum_{l=1}^p P_{l,w}^{(0)}(X_{ij}^{(l)})$ and $R_{ij} = Z_{ij}^*(Y_{ij}^* - Z_{ij}^{*\tau}\beta)$.

THEOREM 3.1.    *Assume that Assumptions 3.1–3.6 hold. Then under* (3.1),

$$(3.7) \qquad \sqrt{mn}[(\hat{\beta} - \beta) - \mu_\beta] \xrightarrow{D} N(0, \Sigma_\beta),$$

*with* $\mu_\beta = (B^{ZZ})^{-1}\mu_B$ *and* $\Sigma_\beta = (B^{ZZ})^{-1}\Sigma_B((B^{ZZ})^{-1})^\tau$, *where* $B^{ZZ} = EZ_{11}^*Z_{11}^{*\tau}$, $\mu_B = E[R_{ij}]$ *and* $\Sigma_B = \sum_{i=-\infty}^\infty \sum_{j=-\infty}^\infty E[(R_{00} - \mu_B)(R_{ij} - \mu_B)^\tau]$.

*Furthermore, when* (1.2) *holds, we have*

$$\mu_\beta = 0, \Sigma_\beta = (B^{ZZ})^{-1}\Sigma_B((B^{ZZ})^{-1})^\tau,$$

*where* $\Sigma_B = \sum_{i=-\infty}^\infty \sum_{j=-\infty}^\infty E[R_{00}R_{ij}^\tau]$, *with* $R_{ij} = Z_{ij}^*\varepsilon_{ij}$ *and* $\varepsilon_{ij} = Y_{ij} - m_0(X_{ij}, Z_{ij}) = Y_{ij} - \mu - Z_{ij}^\tau\beta - g(X_{ij})$.

REMARK 3.2.    Note that

$$\sum_{l=1}^p P_{l,w}^{(0)}(X_{ij}^{(l)}) - \beta^\tau \sum_{l=1}^p P_{l,w}^Z(X_{ij}^{(l)}) = \sum_{l=1}^p (P_{l,w}^{(0)}(X_{ij}^{(l)}) - \beta^\tau P_{l,w}^Z(X_{ij}^{(l)}))$$

$$= \sum_{l=1}^p P_{l,w}(X_{ij}^{(l)}, \beta) \equiv g_a(X_{ij}, \beta).$$

Therefore, $Y_{ij}^* - Z_{ij}^{*\tau}\beta = \varepsilon_{ij} + g(X_{ij}) - g_a(X_{ij}, \beta)$, where $g(X_{ij}) - g_a(X_{ij}, \beta)$ is the residual due to the additive approximation. When (1.2) holds, it means that $g(X_{ij})$ in (1.1) has the expression $g(X_{ij}) = \sum_{l=1}^p g_l(X_{ij}^{(l)}) = \sum_{l=1}^p P_{l,w}(X_{ij}^{(l)}, \beta) = g_a(X_{ij}, \beta)$ and $H(X_{ij}) = \sum_{l=1}^p P_{l,w}^Z(X_{ij}^{(l)})$, and hence, $Y_{ij}^* - Z_{ij}^{*\tau}\beta = \varepsilon_{ij}$. As $\beta$ minimizes $L(\beta) = E[Y_{ij} - m_0(X_{ij}, Z_{ij})]^2$, we have $L'(\beta) = 0$ and $E[\epsilon_{ij}Z_{ij}^*] = E[\epsilon_{ij}(Z_{ij} - E[Z_{ij}|X_{ij}])] = 0$ when (1.2) holds. This implies $E[R_{ij}] = 0$ and, hence, $\mu_\beta = 0$ in (3.7) when the marginal integration estimation procedure is employed for the additive form of $g(\cdot)$.

In both theory and practice, we need to test whether $H_0 : \beta = \beta_0$ holds for a given $\beta_0$. The case where $\beta_0 \equiv 0$ is an important one. Before we state the next result, one needs to introduce some notation. Let

$$\widehat{B}^{ZZ} = \frac{1}{mn} \sum_{i=1}^{m} \sum_{j=1}^{n} \widehat{Z}_{ij}^*(\widehat{Z}_{ij}^*)^\tau, \qquad \widehat{Z}_{ij}^* = \tilde{Z}_{ij} - \sum_{l=1}^{p} \hat{P}_{l,w}^Z(X_{ij}^{(l)}),$$

$$\hat{\mu}_B = \frac{1}{mn} \sum_{i=1}^{m} \sum_{j=1}^{n} \widehat{R}_{ij}, \qquad \widehat{R}_{ij} = \widehat{Z}_{ij}^*(\widehat{Y}_{ij}^* - (\widehat{Z}_{ij}^*)^\tau \hat{\beta}),$$

$$\hat{\mu}_\beta = (\widehat{B}^{ZZ})^{-1} \hat{\mu}_B, \qquad \hat{\Sigma}_\beta = (\widehat{B}^{ZZ})^{-1} \hat{\Sigma}_B ((\widehat{B}^{ZZ})^{-1})^\tau,$$

in which $\hat{\Sigma}_B$ is a consistent estimator of $\Sigma_B$, defined simply by

$$\hat{\Sigma}_B = \sum_{i=-M_m}^{M_m} \sum_{j=-N_n}^{N_n} \hat{\gamma}_{ij},$$

$$\hat{\gamma}_{ij} = \begin{cases} \dfrac{1}{mn} \displaystyle\sum_{u=1}^{m-i} \sum_{v=1}^{n-j} (\widehat{R}_{uv} - \hat{\mu}_B)(\widehat{R}_{u+i,v+j} - \hat{\mu}_B)^\tau, & \text{if (1.1) holds,} \\[2ex] \dfrac{1}{mn} \displaystyle\sum_{u=1}^{m-i} \sum_{v=1}^{n-j} \widehat{R}_{uv} \widehat{R}_{u+i,v+j}^\tau, & \text{if (1.2) holds,} \end{cases}$$

where $M_m \to \infty$, $N_n \to \infty$, $M_m/m \to 0$ and $N_n/n \to 0$ as $m \to \infty$ and $n \to \infty$. It can be shown that both $\hat{\mu}_\beta$ and $\hat{\Sigma}_\beta$ are consistent estimators of $\mu_\beta$ and $\Sigma_\beta$, respectively.

We are now in the position to state a corollary of Theorem 3.1 that can be used to test hypotheses about $\beta$.

COROLLARY 3.1. *Assume that the conditions of Theorem 3.1 hold. Then under (3.1),*

(3.8) $$\hat{\Sigma}_\beta^{-1/2} \sqrt{mn}[(\hat{\beta} - \beta) - \hat{\mu}_\beta] \xrightarrow{D} N(0, I_q)$$

*and*

(3.9) $$mn[(\hat{\beta} - \beta) - \hat{\mu}_\beta]^\tau \hat{\Sigma}_\beta^{-1}[(\hat{\beta} - \beta) - \hat{\mu}_\beta] \xrightarrow{D} \chi_q^2.$$

*Furthermore, when (1.2) holds, we have, under (3.1),*

(3.10) $$\hat{\Sigma}_\beta^{-1/2} \sqrt{mn}(\hat{\beta} - \beta) \xrightarrow{D} N(0, I_q)$$

*and*

(3.11) $$(\sqrt{mn}(\hat{\beta} - \beta))^\tau \hat{\Sigma}_\beta^{-1} (\sqrt{mn}(\hat{\beta} - \beta)) \xrightarrow{D} \chi_q^2.$$

The proof of Theorem 3.1 is relegated to the Appendix, while the proof of Corollary 3.1 is straightforward and therefore omitted.

REMARK 3.3.    Theorem 3.1 implies that there is a big difference between the asymptotic variances in the spatial case and in the time series case. The difference is mainly because the time series is unilateral, while the spatial process is not. Let us consider the simplest case of a line process with $p = q = 1$. In the corresponding time series case where $Y_t = \beta Y_{t-1} + g(Y_{t-2}) + e_t$, $e_t$ is usually assumed to be independent of the past information $\{Y_s, \ s < t\}$; then with $Z_t = Y_{t-1}$ and $X_t = Y_{t-2}$, $\varepsilon_t = Y_t - E(Y_t | X_t, Z_t) = e_t$, therefore $R_t = Z_t^* \varepsilon_t = Z_t^* e_t$ (with $Z_t^*$ defined analogously to $Z_{ij}^*$) is a martingale process with $E[R_0 R_t] = 0$ for $t \neq 0$, which leads to $\Sigma_B = E[R_0^2]$. However, in the bilateral case on the line with the index taking values in $\mathbb{Z}^1$ where $Y_t = \beta Y_{t-1} + g(Y_{t+1}) + e_t$, $e_t$ cannot be assumed to be independent of $(Y_{t-1}, Y_{t+1})$ even when $e_t$ itself is an i.i.d. normal process and $g$ is linear, since under some suitable conditions, as shown in [36], the linear stationary solution may be of the form $Y_t = \sum_{j=-\infty}^{\infty} a_j e_{t-j}$, with all $a_j$ nonzero. Then with $Z_t = Y_{t-1}$ and $X_t = Y_{t+1}$, $\varepsilon_t = Y_t - E(Y_t | X_t, Z_t) \neq e_t$, and usually $E[R_0 R_t] \neq 0$ for $t \neq 0$, which leads to $\Sigma_B \neq E[R_0^2]$.

Next we state the result for the nonparametric component.

THEOREM 3.2.    Assume that Assumptions 3.1–3.6 hold. Then under (3.1), for $x_k \in [-L_k, L_k]$,

$$(3.12) \qquad \sqrt{mnb_k}(\widehat{\widehat{P}}_{k,w}(x_k) - P_{k,w}(x_k) - \text{bias}_{1k}) \xrightarrow{D} N(0, \text{var}_{1k}),$$

where

$$\text{bias}_{1k} = \frac{1}{2} b_k^2 \mu_2(K) \int w_{(-k)}(x^{(-k)}) f_{(-k)}(x^{(-k)}) \frac{\partial^2 g(x, \beta)}{\partial x_k^2} \, dx^{(-k)}$$

and

$$\text{var}_{1k} = J \int V(x, \beta) \frac{[w_{(-k)}(x^{(-k)}) f_{(-k)}(x^{(-k)})]^2}{f(x)} \, dx^{(-k)},$$

with $J = \int K^2(u) \, du$, $\mu_2(K) = \int u^2 K(u) \, du$, $g(x, \beta) = E[(Y_{ij} - \mu - Z_{ij}^\top \beta) | X_{ij} = x]$ and $V(x, \beta) = E[(Y_{ij} - \mu - Z_{ij}^\top \beta - g(x, \beta))^2 | X_{ij} = x]$.

Furthermore, assume that the additive form (1.2) holds and that $E[w_{(-k)}(X_{ij}^{(-k)})] = 1$. Then under (3.1),

$$(3.13) \qquad \sqrt{mnb_k}(\hat{g}_k(x_k) - g_k(x_k) - \text{bias}_{2k}) \xrightarrow{D} N(0, \text{var}_{2k}),$$

*where*

$$\mathrm{bias}_{2k} = \frac{1}{2} b_k^2 \mu_2(K) \frac{\partial^2 g_k(x_k)}{\partial x_k^2}$$

*and*

$$\mathrm{var}_{2k} = J \int V(x,\beta) \frac{[w_{(-k)}(x^{(-k)}) f_{(-k)}(x^{(-k)})]^2}{f(x)} \, dx^{(-k)},$$

*with* $V(x,\beta) = E[(Y_{ij} - \mu - Z_{ij}^\tau \beta - \sum_{k=1}^p g_k(x_k))^2 | X_{ij} = x]$.

The proof of Theorem 3.2 is relegated to the Appendix. We finally state the corresponding results of Theorems 3.1 and 3.2 under Assumptions 3.1–3.3 and 3.4′–3.6′ in Theorem 3.3 below. Its proof is omitted.

THEOREM 3.3. (i) *Assume that Assumptions 3.1–3.3 and 3.4′–3.6′ hold. Then under* (3.1), *the conclusions of Theorem 3.1 hold.*

(ii) *Assume that Assumptions 3.1–3.3 and 3.4′–3.6′ hold. Then under* (3.1), *for* $x_k \in [-L_k, L_k]$,

$$(3.14) \quad \sqrt{mnb_k} (\widehat{P}_{k,w}(x_k) - P_{k,w}(x_k) - \mathrm{bias}_{1k}(I)) \xrightarrow{D} N(0, \mathrm{var}_{1k}(I)),$$

*where*

$$\mathrm{bias}_{1k}(I) = \frac{1}{2} b_k^I \mu_I(K) \int w_{(-k)}(x^{(-k)}) f_{(-k)}(x^{(-k)}) \frac{\partial^I g(x,\beta)}{\partial x_k^I} \, dx^{(-k)}$$

*and*

$$\mathrm{var}_{1k}(I) = J \int V(x,\beta) \frac{[w_{(-k)}(x^{(-k)}) f_{(-k)}(x^{(-k)})]^2}{f(x)} \, dx^{(-k)},$$

*with* $g(x,\beta) = E[(Y_{ij} - \mu - Z_{ij}^\tau \beta)|X_{ij} = x]$, $V(x,\beta) = E[(Y_{ij} - \mu - Z_{ij}^\tau \beta - g(x,\beta))^2 | X_{ij} = x]$, $J = \int K^2(u) \, du$ *and* $\mu_I(K) = \int u^I K(u) \, du$.

*Furthermore, let the additive form* (1.2) *hold and* $E[w_{(-k)}(X_{ij}^{(-k)})] = 1$. *Then under* (3.1),

$$(3.15) \quad \sqrt{mnb_k} (\hat{g}_k(x_k) - g_k(x_k) - \mathrm{bias}_{2k}(I)) \xrightarrow{D} N(0, \mathrm{var}_{2k}(I)),$$

*where*

$$\mathrm{bias}_{2k}(I) = \frac{1}{2} b_k^I \mu_I(K) \frac{\partial^I g_k(x_k)}{\partial x_k^I}$$

*and*

$$\mathrm{var}_{2k}(I) = J \int V(x,\beta) \frac{[w_{(-k)}(x^{(-k)}) f_{(-k)}(x^{(-k)})]^2}{f(x)} \, dx^{(-k)},$$

*with* $V(x,\beta) = E[(Y_{ij} - \mu - Z_{ij}^\tau \beta - \sum_{k=1}^p g_k(x_k))^2 | X_{ij} = x]$.

**4. An illustrative example with simulation.**  In this section we consider an application to the wheat data set of Mercer and Hall [26] as an illustration of the theory and methodology established in this paper. This data set has been analyzed by several investigators including Whittle [36] and Besag [1]; see also [25] on the analysis from the spectral perspective. It involves 500 wheat plots, each 11 ft by 10.82 ft, arranged in a 20×25 rectangle, plot totals constituting the observations. Two measurements, grain yield and straw yield, were made on each plot. Whittle [36] analyzed the grain yields, fitting various stationary unconditional normal autoregressions. Besag [1] analyzed the same data set, but on the basis of the homogenous first- and second-order auto-normal schemes [see (5.5) and (5.6) in [1], page 206], and found that the first-order auto-normal scheme appears satisfactory ([1], page 221). This model has the conditional mean of $Y_{ij}$, given all other site values, equal to

$$(4.1) \qquad \gamma_0 + \gamma_1(Y_{i-1,j} + Y_{i+1,j}) + \gamma_2(Y_{i,j-1} + Y_{i,j+1}),$$

where we use $Y_{ij}$ to denote the grain yield, and $\gamma_0$, $\gamma_1$ and $\gamma_2$ are unknown parameters. For more details, the reader is referred to the above references.

As a first step, we are concerned with whether or not the first-order scheme is linear as in (4.1) or partially linear as in (1.2). This suggests considering the additive first-order scheme

$$(4.2) \qquad \mu + g_1(X_{ij}^{(1)}) + g_2(X_{ij}^{(2)}),$$

where $X_{ij}^{(1)} = Y_{i-1,j} + Y_{i+1,j}$, $X_{ij}^{(2)} = Y_{i,j-1} + Y_{i,j+1}$, $\mu$ is an unknown parameter and $g_1(\cdot)$ and $g_2(\cdot)$ are two unknown functions on $\mathbb{R}^1$. If the Besag scheme is correct, both (1.1) and (1.2) hold and are linear, and one can model (4.2) as a special case of model (1.2) with $\beta = 0$.

Next, we apply the approach established in this paper to estimate $g_1$ and $g_2$. In doing so, the two bandwidths, $b_1 = 0.6$ and $b_2 = 0.7$, were selected using a cross-validation selection procedure for the case of $p = 2$. The resulting estimated functions of $g_1(\cdot)$ and $g_2(\cdot)$ are depicted in Figure 1(a) and (b) with solid lines, respectively, where the additive modeling, based on the modified backfitting algorithm proposed by Mammen, Linton and Nielsen [24] in the i.i.d. case and developed by Lu et al. [23] for the spatial process, is also plotted with dotted lines. We need to point out that, in an asymptotic analysis of such a two-dimensional model, two bandwidths tending to zero at different rates have to be used for each component, thus, we will need to use four bandwidths altogether. But in a finite sample situation like ours, we think that it may be better to rely on cross-validation. This technique is certainly used in the nonspatial situation too, even in cases where an optimal asymptotic formula exists.

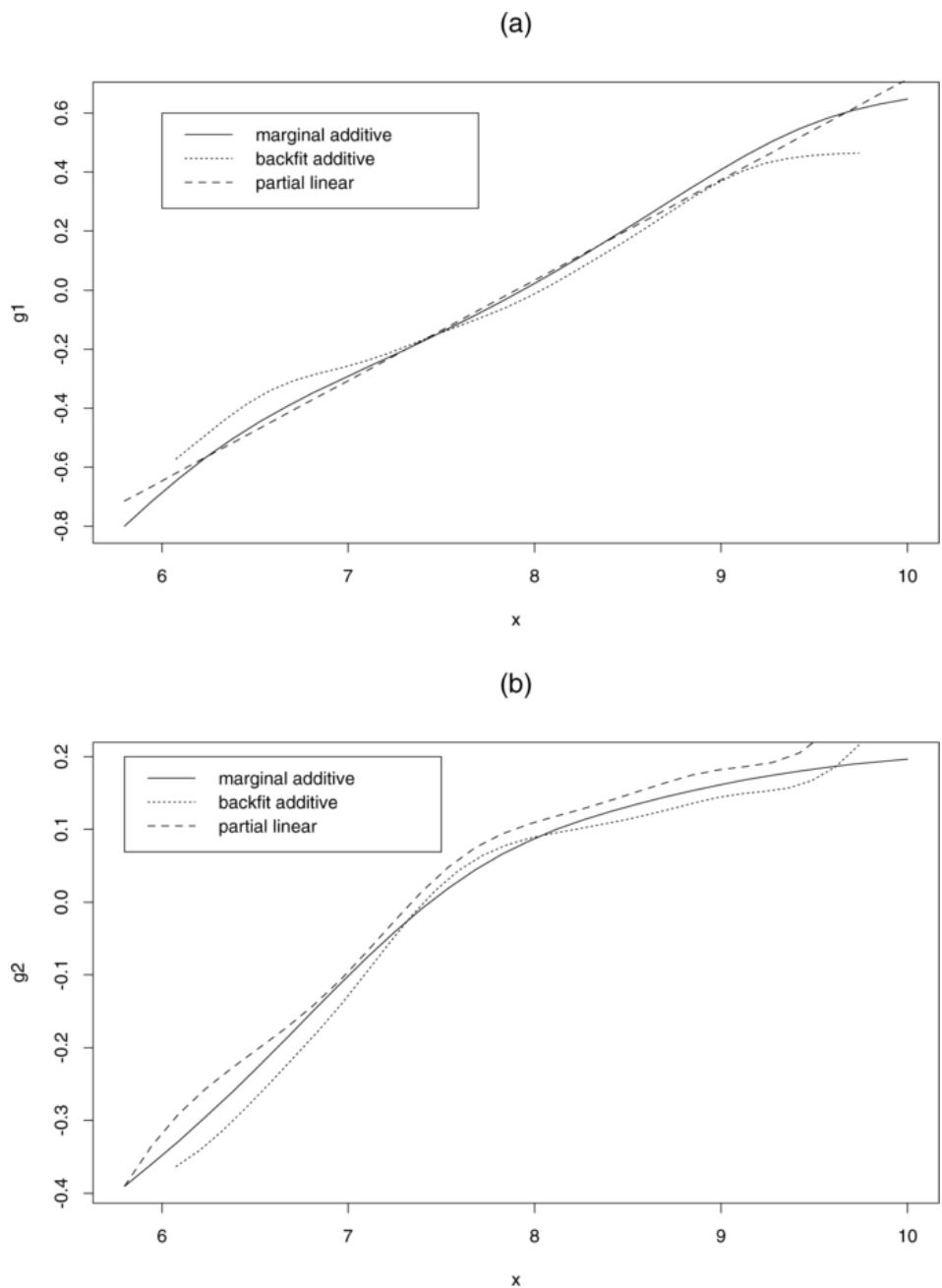

FIG. 1. *Estimated functions of semi-parametric first-order schemes:* (a) $g_1(x)$, (b) $g_2(x)$. *Here the solid and the dotted lines are for the estimates of the additive first-order scheme based on the marginal integration developed in this paper and the modified backfitting in* [24] *and* [23], *respectively; the dashed line is for the estimates of the partially linear first-order scheme based on the approach developed in this paper.*

The pictures of the additive first-order scheme indicate that the estimated function of $g_1(\cdot)$ appears to be linear as in [1], while the estimated function of $g_2(\cdot)$ seems to be nonlinear. This suggests using a partially linear spatial autoregression of the form

$$(4.3) \qquad \qquad \beta_0 + \beta_1 X_{ij}^{(1)} + g_2(X_{ij}^{(2)}).$$

For this case, we view model (4.3) as a special case of model (1.2) with $\mu = \beta_0$, $\beta = \beta_1$, $Z_{ij} = X_{ij}^{(1)}$, $X_{ij} = X_{ij}^{(2)}$ and $g(\cdot) = g_2(\cdot)$. Based on the bandwidth of 0.4 selected using a cross-validation selection procedure, the resulting estimates were $\hat{\beta}_0 = 1.311$, $\hat{\beta}_1 = 0.335$ and $\hat{g}_2(\cdot)$, which are also plotted in Figure 1(a) and (b) with dashed lines, respectively.

We find that our estimate of $\beta_1$ based on the partially linear first-order scheme is almost the same as Besag's first-order auto-normal schemes, which are tabulated in Table 1 below. The estimate of $g_2(\cdot)$ based on the partially linear first-order scheme, similarly to that given in Figure 1(b) based on both the marginal integration and the backfitting of the additive first-order scheme, indicates nonlinearity with a change point around $x = 7.8$.

One may wonder whether the apparent nonlinearity in $g_2$ could arise from random variation even if $g_2$ is linear. The similarity of the two estimates using different techniques is reassuring, but we also did some simulations with samples from the auto-normal first-order scheme with conditional mean (4.1) with $\gamma_0 = 0.16$, $\gamma_1 = 0.34$, $\gamma_2 = 0.14$ and with constant conditional variance $\sigma^2 = 0.11$, where the values of the parameters were chosen to be close to the estimated values of the auto-normal first-order scheme for the grain yields data given by Besag's [1] coding method. The sample size in the simulation is the same as that of the grain yields data, that is, $m = 20$ and $n = 25$. We repeated the simulation 100 times. For each simulated realization, our partially linear first-order scheme of (4.3) was estimated by the approach developed in this paper with the bandwidth of 0.4 (the same as that used for the grain yields data in the above). The boxplots of the 100 simulations for the nonparametric component $g_2(\cdot)$ are depicted in Figure 2. A six-number summary for $\hat{\beta}_1$ is given in Table 2.

TABLE 1
*Estimates of different first-order conditional autoregression schemes for Mercer and Hall's data*

| Scheme | Regressor: $X_{ij}^{(1)}$ | Regressor: $X_{ij}^{(2)}$ | Variance of residuals |
|---|---|---|---|
| Partially linear | $\hat{\beta}_1 = 0.335$ | $\hat{g}_2(\cdot)$: Figure 1(b) | 0.1081 |
| Auto-normal ([1], Table 8) | $\hat{\gamma}_1 = 0.343$ | $\hat{\gamma}_2 = 0.147$ | 0.1099 |
| Auto-normal ([1], Table 10) | $\hat{\gamma}_1 = 0.350$ | $\hat{\gamma}_2 = 0.131$ | 0.1100 |

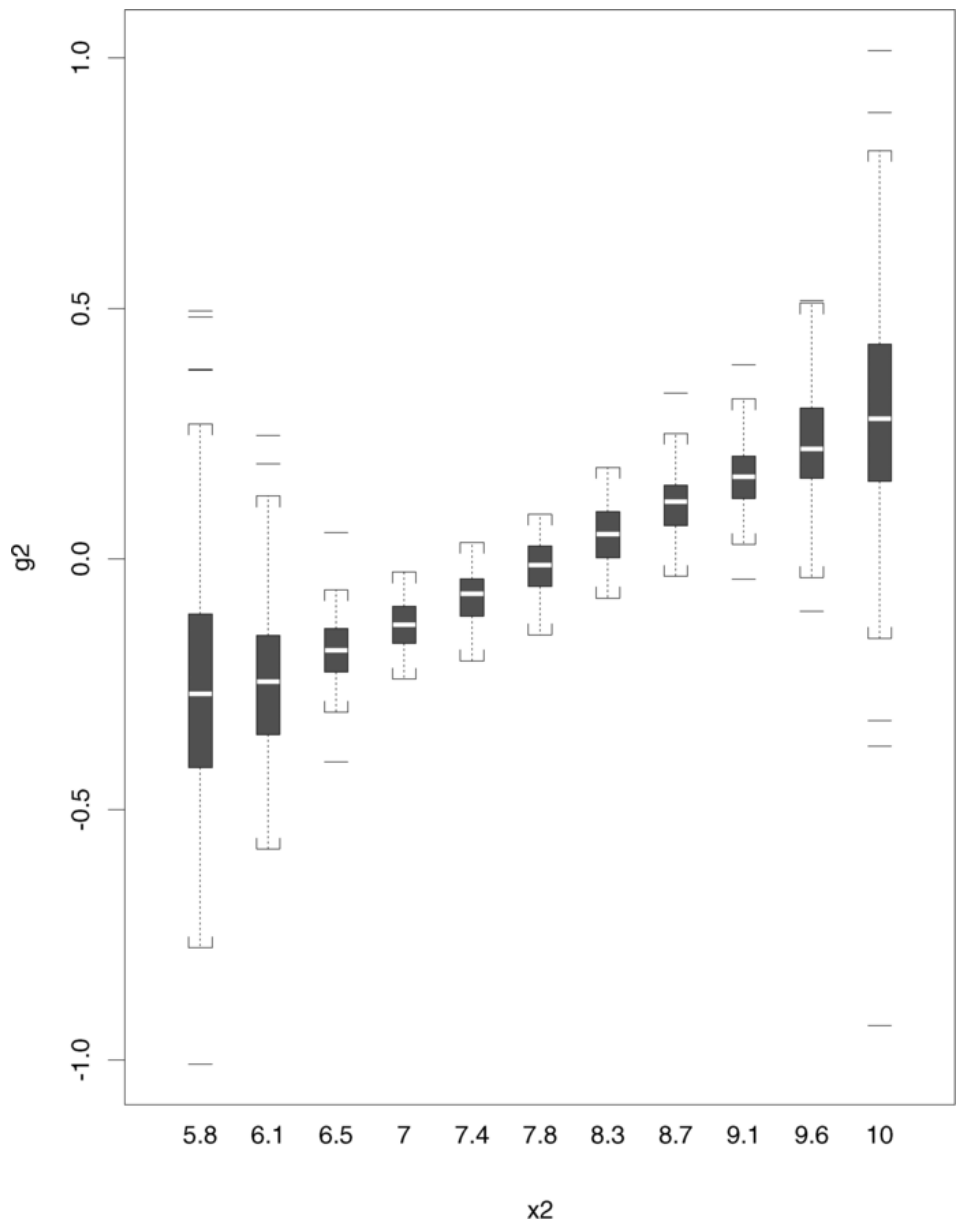

Fig. 2. *Boxplots of the estimated partial linear first-order scheme for the 100 simulations of the auto-normal first-order model for the nonparametric component $g_2(x)$. The sample size is $m = 20$ and $n = 25$.*

It is clear that the estimate for $\beta_1$ is quite stable with median almost equal to the actual parameter $\hat{\beta}_1 = 0.34$, and the estimate for $g_2$ also looks quite linear with small errors around $x = 7.8$. The simulation results show that it is unlikely that the estimated nonlinearity in $g_2$ for the grain yields data in Figure 1(b) should be caused by random variations with the true model being linear. In fact, the accuracy of our estimates is quite high around $x = 7.8$, since the samples of the grain yields are quite dense there (see Figure 3).

Table 2
*A six-number summary for $\hat{\beta}_1$*

| Min. | 1st Qu. | Median | Mean | 3rd Qu. | Max. |
|------|---------|--------|------|---------|------|
| 0.2313 | 0.3129 | 0.3405 | 0.3387 | 0.3684 | 0.4182 |

Table 1 reports the variance of the residuals of the partially linear first-order scheme, as well as of Besag's auto-normal schemes. By contrast, the partially linear first-order scheme gives some improvement over the auto-normal schemes, but perhaps surprisingly small in view of the rather pronounced nonlinearity of Figure 1. In an attempt to understand this, we also calculated the variances of the estimated components and the variance of $Y_{ij}$ over $\{(i, j) : 2 \le i \le 19, 2 \le j \le 24\}$, reported in Table 3. By combining Table 3 with Table 1, we can see the following: (a) clearly, for the partially linear first-order scheme, as well as Besag's auto-normal schemes, the variances of the residuals (in Table 1) are quite large, all about half of the variance of $Y_{ij}$ (given in Table 3); (b) the variances of the first component, $\mathrm{Var}\{g_1(X_{ij}^{(1)})\}$, are much larger (6 times) than those of the second component, $\mathrm{Var}\{g_2(X_{ij}^{(2)})\}$, and therefore, the first components in the fitted conditional means play a key role, while the impact of the second components is smaller; and (c) if we are only concerned with the estimate of the second component $g_2$, then the improvement of the partially linear first-order scheme over the auto-normal schemes is clear if measured in terms of the relative increase of the variance: $(0.0114 - 0.0102)/0.0102 \times 100\% = 11.76\%$ and $(0.0114 - 0.0081)/0.0081 \times 100\% = 40.74\%$ (cf. Table 3). These facts serve at least as tentative explanations of the slightly contradictory messages of

Table 3
*Variances of components of different first-order conditional autoregression schemes for Mercer and Hall's data*

| Scheme | $\mathrm{Var}(Y_{ij})$ | $\mathrm{Var}\{g_1(X_{ij}^{(1)})\}$ | $\mathrm{Var}\{g_2(X_{ij}^{(2)})\}$ |
|--------|------------------------|-------------------------------------|-------------------------------------|
| Partially linear | 0.205 | 0.0661 | 0.0114 |
| Auto-normal ([1], Table 8) | 0.205 | 0.0693 | 0.0102 |
| Auto-normal ([1], Table 10) | 0.205 | 0.0722 | 0.0081 |

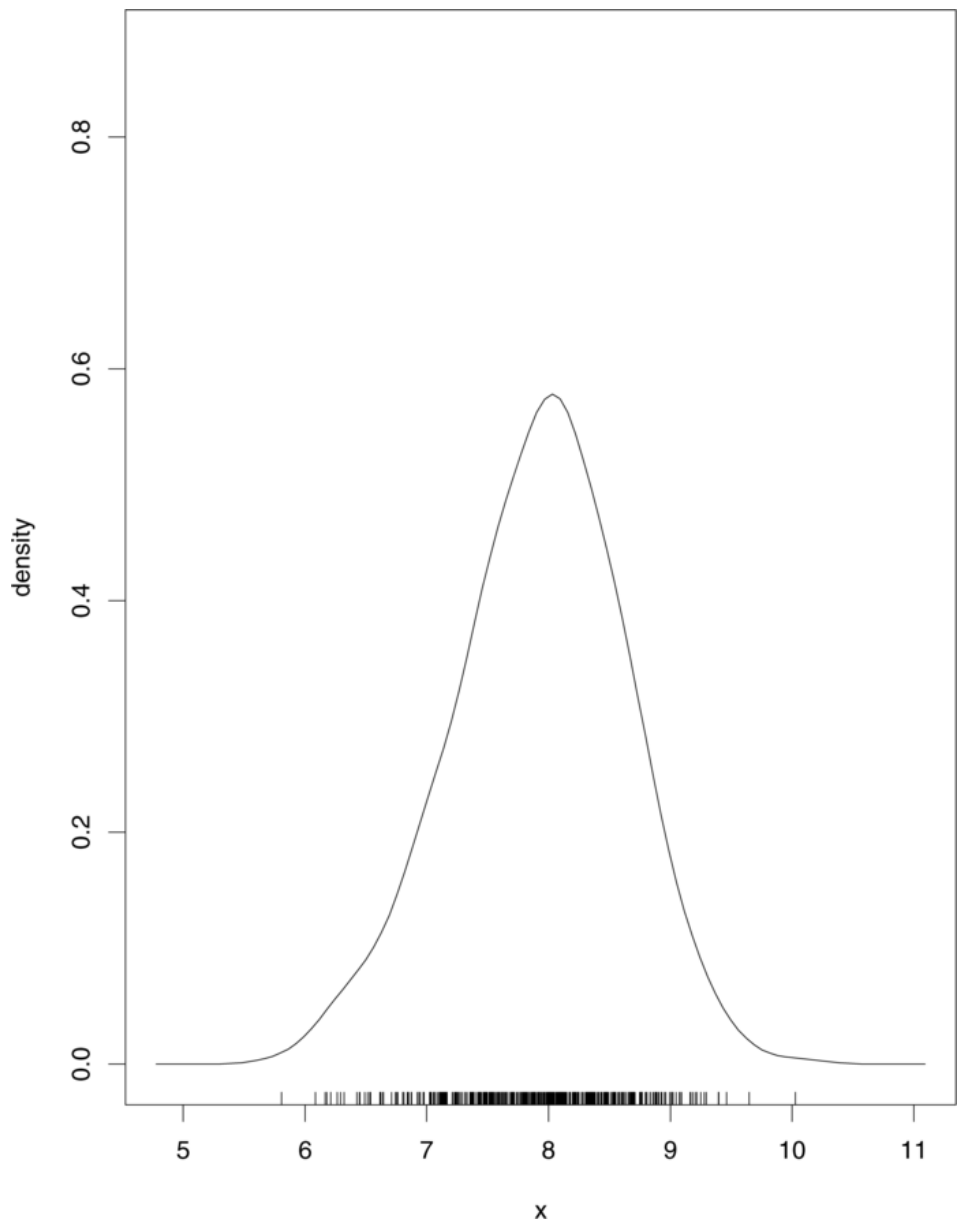

Fig. 3. *The estimated kernel density of $X_{ij}^{(2)}$ defined in (4.3) for the grain yields data.*

Figure 1 and Table 1. The partially linear scheme provides an alternative choice of fitting and conveys more information on the data. A referee suggested that the apparent nonlinearity may be due to an inhomogeneity in the data (cf. [25]). This is a possibility that cannot be ruled out. Also, for time series it is sometimes difficult to distinguish between nonlinearity and nonstationarity.

**5. Conclusion and future studies.** This paper uses a semiparametric additive technique to estimate conditional means of spatial data. The key idea is that the semiparametric technique is employed as an approximation to the true conditional mean function of the spatial data. The asymptotic properties of the resulting estimates are given in Theorems 3.1–3.3. The results of this paper can serve as a starting point for research in a number of directions, including problems related to the estimation of the conditional variance function of a set of spatial data.

In Section 4 our empirical studies show that the estimated form of $g_2(\cdot)$ is nonlinear. To further support such nonlinearity, one may need to establish a formal test. In general, we may consider testing for linearity in the nonparametric components $g_l(\cdot)$ involved in model (1.2).

In the time series case, such test procedures for linearity have been studied extensively during the last ten years. Details may be found in [10]. In the spatial case, Lu et al. [23] propose a bootstrap test and then discuss its implementation. To the best of our knowledge, there is no asymptotic theory available for such a test, and the theoretical problems are very challenging.

To test $H_0 : g_k(X_{ij}^{(k)}) = X_{ij}^{(k)} \gamma_k$, where $\{\gamma_k\}$ is an unknown parameter for each given $k$, our experience with the nonspatial case suggests using a kernel-based test statistic of the form

$$L_k = \sum_{i_1=1}^{m} \sum_{j_1=1}^{n} \sum_{i_2=1,\neq i_1}^{m} \sum_{j_2=1,\neq j_1}^{n} K_{i_1 j_1}(X_{i_2 j_2}, b) \hat{\epsilon}_{i_1 j_1}^{(k)} \hat{\epsilon}_{i_2 j_2}^{(k)},$$

where $K_{i_1 j_1}(X_{i_2 j_2}, b) = \prod_{l=1}^{p} K(\frac{X_{i_1 j_1}^{(l)} - X_{i_2 j_2}^{(l)}}{b_l})$, as defined at the beginning of Section 2, and $\hat{\epsilon}_{ij}^{(k)} = Y_{ij} - \hat{\mu} - Z_{ij}^{\tau} \hat{\beta} - X_{ij}^{(k)} \hat{\gamma}_k - \sum_{l=1,\neq k} \hat{g}_l(X_{ij}^{(l)})$, in which $\hat{\mu}$, $\hat{\beta}$, $\hat{\gamma}_k$ and $\hat{g}_l(\cdot)$ are the corresponding estimators of $\mu$, $\beta$, $\gamma_k$ and $g_l(\cdot)$. These estimators may be defined similarly as in Section 2.

Our experience and knowledge with the nonspatial case would suggest that the normalized version of $L_k$ should have an asymptotically normal distribution under $H_0$, although we have not been able to rigorously prove such a result. This issue and other related issues, for example, a test for isotropy, are left for future research.

## APPENDIX: PROOFS OF THEOREMS 3.1 AND 3.2

Throughout the rest of the paper, the letter $C$ is used to denote constants whose values are unimportant and may vary from line to line. All limits are taken as $(m, n) \to \infty$ in sense of (3.1) unless stated otherwise.

**A.1. Technical lemmas.** In the proofs we need to repeatedly use the following cross term inequality and uniform-consistency lemmas.

Let $f_{(-k)}(\cdot)$ and $f(\cdot)$ be the probability density functions of $X_{ij}^{(-k)}$ and $X_{ij}$, respectively. For $k = 1, 2, \ldots, p$ and $s = 1, 2, \ldots, q$, let

$$d_{ijk}(x_k) = f(X_{ij}^{(-k)}, x_k)^{-1} w(X_{ij}^{(-k)}) f_{(-k)}(X_{ij}^{(-k)}),$$

$$\epsilon_{ij}^{(s)} = Z_{ij}^{(s)} - E[Z_{ij}^{(s)} | X_{ij}], \qquad \Delta_{ij}(x_k) = K\left(\frac{X_{ij}^{(k)} - x_k}{b_k}\right) d_{ijk}(x_k) \epsilon_{ij}^{(s)}.$$

LEMMA A.1. (i) *Let Assumptions* 3.1–3.6 *hold. Then under* (3.1),

$$\frac{1}{\sqrt{mnb_k}} \sum_{i=1}^{m} \sum_{j=1}^{n} \Delta_{ij}(x_k) \xrightarrow{D} N(0, \mathrm{var}_{1k}^{(s)}),$$

*where*

$$\mathrm{var}_{1k}^{(s)} = J \int V^{(s)}(x) \frac{[w_{(-k)}(x^{(-k)}) f_{(-k)}(x^{(-k)})]^2}{f(x)} \, dx^{(-k)},$$

*in which* $J = \int K^2(u) \, du$, $V^{(s)}(x) = E((Z_{ij}^{(s)} - \mu_Z^{(s)} - H^{(s)}(x))^2 | X_{ij} = x)$ *and* $x^{(-k)}$ *is the* $(p-1)$-*dimensional vector obtained from* $x$ *with the* $k$th *component,* $x_k$, *deleted.*

(ii) *Let Assumptions* 3.1–3.6 *hold. For any* $(m, n) \in \mathbb{Z}^2$, *define two sequences of positive integers* $c_1 = c_{1mn}$ *and* $c_2 = c_{2mn}$ *such that* $1 < c_1 < m$ *and* $1 < c_2 < n$. *For any* $x_k$, *let*

$$\tilde{J}(x_k) = \sum_{i=1}^{m} \sum_{j=1}^{n} \sum_{\substack{i'=1 \\ i' \neq i}}^{m} \sum_{\substack{j'=1 \\ j' \neq j}}^{n} E[\Delta_{ij}(x_k) \Delta_{i'j'}(x_k)], \tag{A.1}$$

$$\tilde{J}_1 = c_1 c_2 mn b_k^{(\lambda r - 2)/(\lambda r + 2) + 1}, \tag{A.2}$$

$$\tilde{J}_2 = C mn b_k^{2/(\lambda r)} \left( \sum_{i=\min(c_1, c_2)}^{\sqrt{m^2 + n^2}} i \varphi(i)^{(\lambda r - 2)/(\lambda r)} \right),$$

*where* $C > 0$ *is a positive constant and* $\lambda > 2$ *and* $r \geq 1$ *are as defined in Assumptions* 3.1 *and* 3.2(ii). *Then for any* $x_k$,

$$|\tilde{J}(x_k)| \leq C[\tilde{J}_1 + \tilde{J}_2]. \tag{A.3}$$

PROOF. The proof of (i) follows similarly from that of Lemma 3.1 of [16], while the proof of (ii) is analogous to that of Lemma 5.2 of [16]. When applying Lemma 3.1, one needs to notice that $E[\epsilon_{ij}^{(s)}] = 0$ and $N = 2$. For the application of Lemma 5.2, we need to take $\delta = \lambda r - 2$, $d = 1$ and $N = 2$ in the lemma. □

Lemma A.2.   Let $(i, j) \in \mathbb{Z}^2$ and $\xi_{ij} = K((X_{ij}^{(1)} - x_1)/b_1, \ldots, (X_{ij}^{(p)} - x_p)/b_p)\theta_{ij}$, where $K(\cdot)$ satisfies Assumption 3.5, and $\theta_{ij} = \theta(X_{ij}, Y_{ij})$, in which $\theta(\cdot, \cdot)$ is a measurable function, satisfy $E[\xi_{ij}] = 0$ and $E[|\theta_{ij}|^{\lambda r}] < \infty$ for a positive integer $r$ and some $\lambda > 2$. In addition, let Assumptions 3.1–3.6 hold. Then there exists a constant $C$ depending on $r$ but depending on neither the distribution of $\xi_{ij}$ nor $b_\pi$ and $(m, n)$ such that

$$(A.4) \qquad E\left[\left(\sum_{i=1}^{m}\sum_{j=1}^{n}\xi_{ij}\right)^{2r}\right] \leq C(mnb_\pi)^r$$

holds for all $p$ sets of bandwidths.

Proof.   The proof of this lemma follows from that of Lemma 6.2 of [11]. □

Lemma A.3.   Let $\{Y_{ij}, X_{ij}\}$ be an $\mathbb{R}^1 \times \mathbb{R}^p$-valued stationary spatial process with the mixing coefficient function $\varphi(\cdot)$ as defined in (3.3). Set $\theta_{ij} = \theta(X_{ij}, Y_{ij})$ and $R(x) = E(\theta_{ij}|X_{ij} = x)$. Assume that $E|\theta_{ij}|^{\lambda r} < \infty$ for some positive integer $r$ and some $\lambda > 2$, and that Assumptions 3.1–3.6 hold. Let $R(x)$ and $f(x)$ be twice differentiable with bounded second-order derivatives on $\mathbb{R}^p$. Then

$$
\begin{aligned}
(A.5) \qquad & \sup_{x \in S_W} \left| (mnb_\pi)^{-1} \sum_{i=1}^{m}\sum_{j=1}^{n}\theta_{ij}\prod_{l=1}^{p}K((X_{ij}^{(l)} - x_l)/b_l) - f(x)R(x) \right| \\
& = O_P\left( (mnb_\pi^{1+2/r})^{-r/(p+2r)} + \sum_{k=1}^{p} b_k^2 \right)
\end{aligned}
$$

holds for all $p$ sets of bandwidths.

Proof.   The lemma follows from Lemma A.3 of [11]. □

Lemma A.4.   Let $U_{m,n}$ be as defined in (2.4). Suppose Assumptions 3.1, 3.2 and 3.4 hold. In addition, if $b_\pi \to 0$ and $mnb_\pi \to \infty$, then uniformly over $x \in S_W$,

$$(A.6) \qquad U_{m,n} \xrightarrow{p} U \equiv f(x)\begin{pmatrix} 1 & \mathbf{0}^\tau \\ \mathbf{0} & \mu_2(K)I_p \end{pmatrix},$$

where $\mathbf{0} = (0, \ldots, 0)^\tau \in \mathbb{R}^p$, $\mu_2(K) = \int u^2 K(u)\,du$, $I_p$ is an identity matrix of order $p$ and $\xrightarrow{P}$ denotes convergence in probability.

Proof.   The proof follows from Lemma A.3. Its details are available from the proof of Lemma 6.4 of [11]. □

**A.2. Proofs of Theorems 3.1 and 3.2.** To prove our main theorems, we will often use the property of the marginal integration estimator, which is to be established here and is of independent interest in some other applications.

Let $H^{(s)}(x) = E[(Z^{(s)} - \mu_Z^{(s)})|X = x]$ be the conditional regression of $Z_{ij}^{(s)} - \mu_Z^{(s)}$ given $X_{ij} = x$, $P_{k,w}^{(s)}(x_k) = E[H^{(s)}(X_{ij}^{(-k)}, x_k)w_{(-k)}(X_{ij}^{(-k)})]$ the weighted marginal integration of $H^{(s)}(x)$, and $H_a^{(s)}(x) = \sum_{k=1}^{p} P_{k,w}^{(s)}(x_k)$ the additive approximation of $H^{(s)}(x)$ based on marginal integrations, for $s = 0, 1, \ldots, q$. The estimates of these functionals were given in Section 2. Let $W(x)$ and $S_W$ be as defined in Lemma A.3. The following lemma is necessary for the proof of the main theorems.

LEMMA A.5. *Suppose Assumptions 3.1–3.5 hold and the bandwidths satisfy $mnb_k^5 = O(1)$, $\sum_{l=1, l \neq k}^{p} b_l^2 = o(b_k^2)$. Then under (3.1),*

$$\text{(A.7)} \qquad \sqrt{mnb_k}(\hat{P}_{k,w}^{(s)}(x_k) - P_{k,w}^{(s)}(x_k) - \text{bias}_{1k}^{(s)}) \xrightarrow{D} N(0, \text{var}_{1k}^{(s)}),$$

*where*

$$\text{bias}_{1k}^{(s)} = \frac{1}{2} b_k^2 \mu_2(K) \int w_{(-k)}(x^{(-k)}) f_{(-k)}(x^{(-k)}) \frac{\partial^2 H^{(s)}(x)}{\partial x_k^2} \, dx^{(-k)},$$

$$\text{var}_{1k}^{(s)} = J \int V^{(s)}(x) \frac{[w_{(-k)}(x^{(-k)}) f_{(-k)}(x^{(-k)})]^2}{f(x)} \, dx^{(-k)},$$

*in which $\mu_2(K) = \int u^2 K(u) \, du$, and the other quantities are as defined in Lemma A.1.*

*Let $H_k^{(s)}(x_k) = E[(Z_{ij}^{(s)} - \mu_Z^{(s)})|X_{ij}^{(k)} = x_k]$. Furthermore, if $H^{(s)}(x) = \sum_{k=1}^{p} H_k^{(s)}(x_k)$ and $E[w_{(-k)}(X_{ij}^{(-k)})] = 1$, then under (3.1),*

$$\text{(A.8)} \qquad \sqrt{mnb_k}(\hat{P}_{k,w}^{(s)}(x_k) - H_k^{(s)}(x_k) - \text{bias}_{2k}^{(s)}) \xrightarrow{D} N(0, \text{var}_{2k}^{(s)}),$$

*where*

$$\text{bias}_{2k}^{(s)} = \frac{1}{2} b_k^2 \mu_2(K) \frac{\partial^2 H_k^{(s)}(x_k)}{\partial x_k^2}$$

*and*

$$\text{var}_{2k}^{(s)} = J \int V^{(s)}(x) \frac{[w_{(-k)}(x^{(-k)}) f_{(-k)}(x^{(-k)})]^2}{f(x)} \, dx^{(-k)},$$

*where $V^{(s)}(x) = E[(Z_{ij}^{(s)} - \mu_Z^{(s)} - \sum_{k=1}^{p} H_k^{(s)}(x_k))^2 | X_{ij} = x]$.*

PROOF.    By the law of large numbers, it is obvious that, for $x_k \in [-L_k, L_k]$,

$$\tilde{P}_{k,w}^{(s)}(x_k) = (mn)^{-1} \sum_{i=1}^{m} \sum_{j=1}^{n} H^{(s)}(X_{ij}^{(-k)}, x_k) w_{(-k)}(X_{ij}^{(-k)})$$

(A.9)

$$= P_{k,w}^{(s)}(x_k) + O_P\left(\frac{1}{\sqrt{mn}}\right).$$

Throughout the rest of the proof, set $\gamma = (1, 0, \ldots, 0)^\tau \in \mathbb{R}^{1+p}$. Note that, by the notation and definitions in Section 2,

$$H_{m,n}^{(s)}(X_{ij}^{(-k)}, x_k) - H^{(s)}(X_{ij}^{(-k)}, x_k)$$

(A.10)

$$= \gamma^\tau U_{m,n}^{-1}(X_{ij}^{(-k)}, x_k) V_{m,n}^{(s)}(X_{ij}^{(-k)}, x_k) - H^{(s)}(X_{ij}^{(-k)}, x_k)$$

$$= \gamma^\tau U_{m,n}^{-1}(X_{ij}^{(-k)}, x_k) B_{m,n}(X_{ij}^{(-k)}, x_k),$$

where $DH^{(s)}(x) = (\partial H^{(s)}(x)/\partial x_1, \ldots, \partial H^{(s)}(x)/\partial x_p)$ with $x = (x^{(-k)}, x_k)$, the symbol $\odot$ is as defined in (2.5) and

$$B_{m,n}(x)$$

(A.11)

$$= \begin{pmatrix} v_{m,n,0}^{(s)}(x) - u_{m,n,00}(x)H^{(s)}(x) - U_{m,n,01}(x)(DH^{(s)}(x) \odot b)^\tau \\ V_{m,n,1}^{(s)}(x) - U_{m,n,10}(x)H^{(s)}(x) - U_{m,n,11}(x)(DH^{(s)}(x) \odot b)^\tau \end{pmatrix}$$

$$\equiv \begin{pmatrix} B_{m,n,0}(x) \\ B_{m,n,1}(x) \end{pmatrix}.$$

Therefore, by the uniform consistency in Lemma A.4, for $x_k \in [-L_k, L_k]$,

$$\hat{P}_{k,w}^{(s)}(x_k) - \tilde{P}_{k,w}^{(s)}(x_k)$$

$$= \gamma^\tau (mn)^{-1} \sum_{i=1}^{m} \sum_{j=1}^{n} U_{m,n}^{-1}(X_{ij}^{(-k)}, x_k) B_{m,n}(X_{ij}^{(-k)}, x_k) w_{(-k)}(X_{ij}^{(-k)})$$

(A.12)

$$= (mn)^{-1} \sum_{i=1}^{m} \sum_{j=1}^{n} f^{-1}(X_{ij}^{(-k)}, x_k) B_{m,n,0}(X_{ij}^{(-k)}, x_k) w_{(-k)}(X_{ij}^{(-k)})$$

$$+ O_P(d_{mn})(mn)^{-1} \sum_{i=1}^{m} \sum_{j=1}^{n} B_{m,n,0}(X_{ij}^{(-k)}, x_k) w_{(-k)}(X_{ij}^{(-k)}),$$

where $d_{mn} = (mnb_\pi^{1+2/r})^{-r/(p+2r)} + \sum_{l=1}^{p} b_l^2$. Note that

$$B_{m,n,0}(x) = (mnb_\pi)^{-1} \sum_{i'=1}^{m} \sum_{j'=1}^{n} \left( \tilde{Z}_{i'j'}^{(s)} - H^{(s)}(x) \right.$$

$$-\sum_{\ell=1}^{p}\frac{\partial H^{(s)}}{\partial x_\ell}(x)(X_{i'j'}^{(\ell)}-x_\ell)\bigg)K_{i'j'}(x,b)$$

(A.13)
$$=(mnb_\pi)^{-1}\sum_{i'=1}^{m}\sum_{j'=1}^{n}\eta_{i'j'}(x)K_{i'j'}(x,b)$$

$$-(\overline{Z}^{(s)}-\mu_Z^{(s)})(mnb_\pi)^{-1}\sum_{i'=1}^{m}\sum_{j'=1}^{n}K_{i'j'}(x,b)$$

$$\equiv B_{m,n,0}^{*}(x^{(-k)},x_k)+B_{m,n,0}^{**}(x^{(-k)},x_k),$$

where $\eta_{i'j'}(x)=Z_{i'j'}^{(s)}-\mu_Z^{(s)}-H^{(s)}(x)-\sum_{l=1}^{p}\frac{\partial H^{(s)}}{\partial x_l}(x)(X_{i'j'}^{l}-x_l)$.

Clearly, the result of $\overline{Z}^{(s)}-\mu_Z^{(s)}=O_P(\frac{1}{\sqrt{mn}})$ together with the uniform consistency in Lemma A.3 leads to

$$B_{m,n,0}^{**}(x^{(-k)},x_k)=O_P\left(\frac{1}{\sqrt{mn}}\right),$$

which holds uniformly with respect to $x=(x^{(-k)},x_k)\in S_W$. Now it follows from (A.12)–(A.13) by exchanging the summations over $(i,j)$ and $(i',j')$ that

$$\hat{P}_{k,w}^{(s)}(x_k)-\tilde{P}_{k,w}^{(s)}(x_k)$$

(A.14)
$$=(mnb_k)^{-1}\sum_{i'=1}^{m}\sum_{j'=1}^{n}K\left(\frac{X_{i'j'}^{(k)}-x_k}{b_k}\right)B_{i'j'}^{(k)}(x_k)$$

$$+O_P(c_{mn})(mnb_k)^{-1}\sum_{i'=1}^{m}\sum_{j'=1}^{n}K\left(\frac{X_{i'j'}^{(k)}-x_k}{b_k}\right)B_{i'j'}^{*(k)}(x_k)$$

$$+O_P\left(\frac{1}{\sqrt{mn}}\right),$$

where $B_{i'j'}^{(k)}(x_k)=\frac{1}{mnb_{(-k)}}\sum_{i=1}^{m}\sum_{j=1}^{n}f^{-1}(X_{ij}^{(-k)},x_k)w_{(-k)}(X_{ij}^{(-k)})\eta_{i'j'}(X_{ij}^{(-k)},x_k)K_{ij,i'j'}^{(-k)}$ and $B_{i'j'}^{*(k)}(x_k)=\frac{1}{mnb_{(-k)}}\sum_{i=1}^{m}\sum_{j=1}^{n}w_{(-k)}(X_{ij}^{(-k)})\eta_{i'j'}(X_{ij}^{(-k)},x_k)K_{ij,i'j'}^{(-k)}$, in which $b_{(-k)}=\prod_{l=1,l\neq k}^{p}b_l$ and $K_{ij,i'j'}^{(-k)}=\prod_{l=1,l\neq k}^{p}K(\frac{X_{ij}^{(l)}-X_{i'j'}^{(l)}}{b_l})$.

Recall $\epsilon_{ij}^{(s)}=Z_{ij}^{(s)}-\mu_Z^{(s)}-H^{(s)}(X_{ij})=Z_{ij}^{(s)}-E(Z_{ij}^{(s)}|X_{ij})$. Note that the properties (compact support) of the kernel function in Assumption 3.5 show that, if $K_{ij,i'j'}^{(-k)}>0$ and $K((X_{i'j'}^{(k)}-x_k)/b_k)>0$ in (A.14), then $|X_{i'j'}^{(l)}-X_{ij}^{(l)}|\leq Cb_l\to0$ for $l\neq k$ and $|X_{i'j'}^{(k)}-x_k|\leq Cb_k\to0$, as $m\to\infty$ and $n\to\infty$. Therefore, if $K_{ij,i'j'}^{(-k)}>0$ and $K((X_{i'j'}^{(k)}-x_k)/b_k)>0$ in (A.14), then by Taylor's expansion (around $X_{ij}$) together with the uniform continuity of second partial derivatives of $g(\cdot)$ in Assumption 3.4,

$$
\begin{aligned}
\eta_{i'j'}(X_{ij}^{(-k)}, x_k) ={}& Z_{i'j'}^{(s)} - \mu_Z^{(s)} - H^{(s)}(X_{ij}^{(-k)}, x_k) \\
& - \sum_{l=1, l\neq k}^{p} \frac{\partial H^{(s)}}{\partial x_l}(X_{ij}^{(-k)}, x_k)(X_{i'j'}^{(\ell)} - X_{ij}^{(l)}) \\
& - \frac{\partial H^{(s)}}{\partial x_k}(X_{ij}^{(-k)}, x_k)(X_{i'j'}^{(k)} - x_k) \\
={}& \epsilon_{i'j'}^{(s)} + \frac{1}{2} \frac{\partial^2 H^{(s)}(X_{ij}^{(-k)}, x_k)}{\partial x_k^2}(X_{i'j'}^{(k)} - x_k)^2 \\
& + \frac{o(1)}{2}\left[\sum_{l,l'=1,\neq k}^{p} b_l b_{l'} + \sum_{l=1,\neq k}^{p} b_l b_k + b_k^2\right] \\
& + \frac{1}{2}\sum_{l,l'=1,\neq k}^{p} \frac{\partial^2 H^{(s)}(X_{ij}^{(-k)}, x_k)}{\partial x_l \partial x_{l'}} O(b_l b_{l'}) \\
& + \sum_{l=1,\neq k}^{p} \frac{\partial^2 H^{(s)}(X_{ij}^{(-k)}, x_k)}{\partial x_l \partial x_k} O(b_l b_k).
\end{aligned}
$$

Then under $K_{ij,i'j'}^{(-k)} > 0$ and $K((X_{i'j'}^{(k)} - x_k)/b_k) > 0$,

$$
\begin{aligned}
B_{i'j'}^{(k)}(x_k) ={}& \epsilon_{i'j'}^{(s)}\{mnb_{(-k)}\}^{-1}\sum_{i=1}^{m}\sum_{j=1}^{n} f^{-1}(X_{ij}^{(-k)}, x_k)w_{(-k)}(X_{ij}^{(-k)})K_{ij,i'j'}^{(-k)} \\
& - \frac{1}{2}(X_{i'j'}^{(k)} - x_k)^2\{mnb_{(-k)}\}^{-1}\sum_{i=1}^{m}\sum_{j=1}^{n} f^{-1}(X_{ij}^{(-k)}, x_k)w_{(-k)}(X_{ij}^{(-k)}) \\
& \hspace{4cm} \times \frac{\partial^2 H^{(s)}(X_{ij}^{(-k)}, x_k)}{\partial x_k^2}K_{ij,i'j'}^{(-k)} \\
& + \frac{1}{2}\sum_{l,l'=1,\neq k}^{p} O(b_l b_{l'})\{mnb_{(-k)}\}^{-1}\sum_{i=1}^{m}\sum_{j=1}^{n} f^{-1}(X_{ij}^{(-k)}, x_k)w_{(-k)}(X_{ij}^{(-k)}) \\
& \hspace{4cm} \times \frac{\partial^2 H^{(s)}(X_{ij}^{(-k)}, x_k)}{\partial x_l \partial x_{l'}}K_{ij,i'j'}^{(-k)} \\
& + \sum_{l=1,\neq k}^{p} O(b_l b_k)\{mnb_{(-k)}\}^{-1}\sum_{i=1}^{m}\sum_{j=1}^{n} f^{-1}(X_{ij}^{(-k)}, x_k)w_{(-k)}(X_{ij}^{(-k)})
\end{aligned}
$$

$$\times \frac{\partial^2 H^{(s)}(X_{ij}^{(-k)}, x_k)}{\partial x_l \, \partial x_k} K_{ij,i'j'}^{(-k)}$$

$$+ \left(\frac{1}{2} \sum_{l,l'=1, \neq k}^{p} b_l b_{l'} + \sum_{l=1, \neq k}^{p} b_l b_k + b_k^2\right) \cdot o(1)$$

$$\times \frac{1}{mnb_{(-k)}} \sum_{i=1}^{m} \sum_{j=1}^{n} f^{-1}(X_{ij}^{(-k)}, x_k) w_{(-k)}(X_{ij}^{(-k)}) K_{ij,i'j'}^{(-k)}.$$

Again, using the uniform consistency in Lemma A.3, we have

$$
\begin{aligned}
B_{i'j'}^{(k)}(x_k) = {} & d_{i'j'k}(x_k) \left[\epsilon_{i'j'}^{(s)} + \frac{1}{2}(X_{i'j'}^{(k)} - x_k)^2 \frac{\partial^2 H^{(s)}(X_{i'j'}^{(-k)}, x_k)}{\partial x_k^2}\right] \\
& + O_P(c_{mn}^{(-k)}) \left[\epsilon_{i'j'}^{(s)} + \frac{1}{2}(X_{i'j'}^{(k)} - x_k)^2\right] \\
& + \frac{1}{2} \sum_{l,l'=1, \neq k}^{p} O(b_l b_{l'}) \left[d_{i'j'k}(x_k) \frac{\partial^2 H^{(s)}(X_{i'j'}^{(-k)}, x_k)}{\partial x_l \, \partial x_{l'}} + O_P(c_{mn}^{(-k)})\right] \\
& + \sum_{l=1, \neq k}^{p} O(b_l b_k) \left[d_{i'j'k}(x_k) \frac{\partial^2 H^{(s)}(X_{i'j'}^{(-k)}, x_k)}{\partial x_l \, \partial x_k} + O_P(c_{mn}^{(-k)})\right] \\
& + \left(\frac{1}{2} \sum_{l,l'=1, \neq k}^{p} o(1) b_l b_{l'} + \sum_{\ell=1, \neq k}^{p} o(1) b_l b_k + o(1) b_k^2\right) \\
& \times [d_{i'j'k}(x_k) + O_P(c_{mn}^{(-k)})],
\end{aligned}
$$

(A.15)

where $d_{ijk}(x_k) = f(X_{ij}^{(-k)}, x_k)^{-1} w_{(-k)}(X_{ij}^{(-k)}) f_{(-k)}(X_{ij}^{(-k)})$.

In addition, denote by

$$d_{ijk}^*(x_k) \equiv w_{(-k)}(X_{ij}^{(-k)}) f_{(-k)}(X_{ij}^{(-k)}) \quad \text{and} \quad K_{b_k}(x_k) \equiv b_k^{-1} K\left(\frac{x_k}{b_k}\right).$$

Then similarly to (A.15),

$$
\begin{aligned}
B_{i'j'k}^{*(k)}(x_k) = {} & d_{i'j'k}^*(x_k) \left[\epsilon_{i'j'}^{(s)} + \frac{1}{2}(X_{i'j'}^{(k)} - x_k)^2 \frac{\partial^2 H^{(s)}(X_{i'j'}^{(-k)}, x_k)}{\partial x_k^2}\right] \\
& + O_P(c_{mn}^{(-k)}) \left[\epsilon_{i'j'}^{(s)} + \frac{1}{2}(X_{i'j'}^{(k)} - x_k)^2\right] \\
& + \frac{1}{2} \sum_{l,l'=1, \neq k}^{p} O(b_l b_{l'}) \left[d_{i'j'k}^*(x_k) \frac{\partial^2 H^{(s)}(X_{i'j'}^{(-k)}, x_k)}{\partial x_l \, \partial x_{l'}} + O_P(c_{mn}^{(-k)})\right]
\end{aligned}
$$

(A.16)

$$+ \sum_{l=1,\neq k}^{p} O(b_l b_k) \left[ d_{i'j'k}^*(x_k) \frac{\partial^2 H^{(s)}(X_{i'j'}^{(-k)}, x_k)}{\partial x_l \, \partial x_k} + O_P(c_{mn}^{(-k)}) \right]$$

$$+ \left( \frac{1}{2} \sum_{l,l'=1,\neq k}^{p} o(1) b_l b_{l'} + \sum_{l=1,\neq k}^{p} o(1) b_l b_k + o(1) b_k^2 \right)$$

$$\times [d_{i'j'k}^*(x_k) + O_P(c_{mn}^{(-k)})].$$

Therefore, by (A.14)–(A.16),

$$\hat{P}_{k,w}^{(s)}(x_k) - \tilde{P}_{k,w}^{(s)}(x_k)$$

(A.17)
$$= T_{mn}^{(k)} + O_P(c_{mn}) T_{mn}^{*(k)} + O_P(1) \sum_{l=1,\neq k}^{p} b_l^2$$

$$+ O_P(1) \sum_{l=1,\neq k}^{p} b_\ell b_k + o_P(1) b_k^2 + O_P(1) \left( \frac{1}{\sqrt{mn}} \right),$$

where

$$T_{mn}^{(k)} = (mnb_k)^{-1} \sum_{i=1}^{m} \sum_{j=1}^{n} K\left( \frac{X_{ij}^{(k)} - x_k}{b_k} \right) d_{ijk}(x_k) \epsilon_{ij}^{(s)}$$

$$+ (mnb_k)^{-1} \sum_{i=1}^{m} \sum_{j=1}^{n} K\left( \frac{X_{ij}^{(k)} - x_k}{b_k} \right) d_{ijk}(x_k)$$

(A.18)
$$\times \left[ \frac{1}{2} (X_{ij}^{(k)} - x_k)^2 \frac{\partial^2 H^{(s)}(X_{ij}^{(-k)}, x_k)}{\partial x_k^2} \right]$$

$$\equiv T_{mn1}^{(k)} + T_{mn2}^{(k)},$$

and $T_{mn}^{*(k)}$ can be expressed similarly to (A.18) with $d_{ijk}(x_k)$ replaced by $d_{ijk}^*(x_k)$.

We next consider $T_{mn1}^{(k)}$ and $T_{mn2}^{(k)}$. Clearly, $E[T_{mn1}^{(k)}] = 0$ since $E(\epsilon_{ij}^{(s)}|X_{ij}) = 0$. We calculate the asymptotic variance of $T_{mn1}^{(k)}$. Note that

(A.19)
$$E[T_{mn1}^{(k)}]^2 = J_1(x_k) + J_2(x_k),$$

where

$$J_1(x_k) = (mnb_k)^{-2} \sum_{i=1}^{m} \sum_{j=1}^{n} E\left[ K^2\left( \frac{X_{ij}^{(k)} - x_k}{b_k} \right) d_{ijk}^2(x_k) (\epsilon_{ij}^{(s)})^2 \right],$$

$$J_2(x_k) = (mnb_k)^{-2} \sum_{i=1}^{m} \sum_{j=1}^{n} \sum_{\substack{i'=1 \\ i' \neq i \text{ or}}}^{m} \sum_{\substack{j'=1 \\ j' \neq j}}^{n} E[\Delta_{ij}(x_k) \Delta_{i'j'}(x_k)],$$

in which $\Delta_{ij}(x_k) = K((X_{ij}^{(k)} - x_k)/b_k) \, d_{ijk}(x_k) \epsilon_{ij}^{(s)}$. A simple calculation implies

$$
\begin{aligned}
\text{(A.20)} \quad J_1(x_k) &= \frac{1}{mnb_k} J E[d_{ijk}^2(x_k) \epsilon_{ij}^2 | X_{ij}^{(k)} = x_k] f_k(x_k)(1 + o(1)) \\
&= \frac{1}{mnb_k}(1 + o(1)) C_k(J, V),
\end{aligned}
$$

where

$$C_k(J, V) = J \int V^{(s)}(x) \frac{[w_{(-k)}(x^{(-k)}) f_{(-k)}(x^{(-k)})]^2}{f(x)} \, dx^{(-k)},$$

in which $J = \int K^2(u) \, du$, $V^{(s)}(x) = E[(\epsilon_{ij}^{(s)})^2 | X_{ij} = x]$, and $f_k(x_k)$ is the density function of $X_{ij}^{(k)}$. To deal with the cross term $J_2(x_k)$, we need to use Lemma 6.1. Under the assumptions of the lemma, it leads to

$$
\begin{aligned}
\text{(A.21)} \quad J_2(x_k) &\leq C(mnb_k)^{-1} \Bigg[ b_k^{(\lambda r - 2)/(\lambda r + 2)} c_1 c_2 \\
&\quad + b_k^{-(\lambda r - 2)/(\lambda r)} \left( \sum_{t=\min\{c_1, c_2\}}^{\infty} t\{\varphi(t)\}^{(\lambda r - 2)/(\lambda r)} \right) \Bigg].
\end{aligned}
$$

Take $c_1 = c_2 = [b_k^{-(\lambda r - 2)/(a\lambda r)}]$, where $[u] \leq u$ denotes the largest integer part of $u$. Then since $a > 2(\lambda r + 2)/\lambda r$ in Assumption 3.3, $\frac{2(\lambda r - 2)}{a\lambda r} < \frac{\lambda r - 2}{\lambda r + 2}$, and it hence follows from (A.21) and Assumption 3.3 that

$$
\begin{aligned}
\text{(A.22)} \quad J_2(x_k) &\leq C(mnb_k)^{-1} \Bigg[ b_k^{(\lambda r - 2)/(\lambda r + 2) - (2(\lambda r - 2))/(a\lambda r)} \\
&\quad + c_1^a \sum_{t=c_1}^{\infty} t\{\varphi(t)\}^{(\lambda r - 2)/(r)} \Bigg] \\
&= o((mnb_k)^{-1}),
\end{aligned}
$$

using $c_1^a \sum_{t=c_1}^{\infty} t\{\varphi(t)\}^{(\lambda r - 2)/(\lambda r)} \leq c_1^a \sum_{t=c_1}^{\infty} t^{2r-1}\{\varphi(t)\}^{(\lambda r - 2)/(\lambda r)} \to 0$ by Assumption 3.3.

Now the asymptotic variance of $T_{mn1}^{(k)}$, using (A.19), (A.20) and (A.22), equals the right-hand side of (A.20), that is,

$$(mnb_k) E[T_{mn1}^{(k)}]^2 \to J \int V^{(s)}(x) \frac{[w_{(-k)}(x^{(-k)}) f_{(-k)}(x^{(-k)})]^2}{f(x)} \, dx^{(-k)}$$

(A.23)
$$\equiv \operatorname{var}_{1k}^{(s)}.$$

Next, we consider the term $T_{mn2}^{(k)}$ in (A.18). From (A.18), together with the property of the kernel function in Assumption 3.5,

$$T_{mn2}^{(k)} = \frac{1}{2} b_k^2 E\left[ d_{ijk}(x_k) \frac{\partial^2 H^{(s)}(X_{ij}^{(-k)}, x_k)}{\partial x_k^2} \Big| X_{ij}^{(k)} = x_k \right] f_k(x_k) \mu_2(K) + O_P(l_{mn}^{(k)}) b_k^2$$

$$= \frac{b_k^2 \mu_2(K) f_k(x_k)}{2} \int w_{(-k)}(x^{(-k)}) \frac{\partial^2 g(x^{(-k)}, x_k)}{\partial x_k^2} \, dx^{(-k)} + o_P(b_k^2)$$

$$\equiv \operatorname{bias}_{1k}^{(s)} + o_P(b_k^2),$$

where $l_{mn}^{(k)} = (mnb_k^{1+2/r})^{-r/(1+2r)} + b_k^2$ and $\mu_2(K) = \int u^2 K(u) \, du$.

Similarly, one can show $T_{mn}^{*(k)} = O_P(1/\sqrt{mnb_k} + b_k^2)$. Based on the conditions, $mnb_k^5 = O(1)$ and $\sum_{\ell=1, \neq k}^p b_\ell^2 = o(b_k^2)$, the remaining terms in (A.17) can be neglected since

$$\sqrt{mnb_k} c_{mn}\left( \frac{1}{\sqrt{mnb_k}} + b_k^2 \right) = (1 + b_k^2 \sqrt{mnb_k}) \left( (mnb_\pi^{1+2/r})^{-r/(p+2r)} + \sum_{l=1}^p b_l^2 \right)$$
$$\to 0,$$

$$\sqrt{mnb_k} \sum_{l=1, \neq k}^p b_l^2 = O(1) \left[ mnb_k \left( \sum_{l=1, \neq k}^p b_l^2 \right)^2 \right]^{1/2} \to 0,$$

$$\sqrt{mnb_k} \sum_{l=1, \neq k}^p b_l b_k = O(1) \left( mnb_k^3 \sum_{l=1, \neq k}^p b_l^2 \right)^{1/2} \to 0$$

and $\sqrt{mnb_k} \frac{1}{\sqrt{mn}} = b_k^{1/2} \to 0$.

Therefore, in view of what we have derived, to complete the proof of (A.8), it suffices to show that $\sqrt{mnb_k} T_{mn1}^{(k)} \xrightarrow{D} N(0, \operatorname{var}_{1k}^{(s)})$, which follows from Lemma A.1(i). $\square$

PROOF OF THEOREM 3.1. We note that

(A.24)
$$\hat{\beta} - \beta = \left( \frac{1}{mn} \sum_{i=1}^m \sum_{j=1}^n \hat{Z}_{ij}^*(\hat{Z}_{ij}^*)^\tau \right)^{-1} \left( \frac{1}{mn} \sum_{i=1}^m \sum_{j=1}^n \hat{Z}_{ij}^*(\hat{Y}_{ij}^* - \hat{Z}_{ij}^* \beta) \right)$$

$$\equiv (B_{mn}^{ZZ})^{-1} B_{mn}^{ZY}.$$

Denote by $H_a^{(s)}(x) \equiv \sum_{l=1}^p P_{l,w}^{(s)}(x_l)$ and $H_a(x) \equiv \sum_{l=1}^p P_{l,w}^Z(x_l)$ the additive approximate versions to $H^{(s)}(x) = E[(Z_{ij}^{(s)} - \mu_Z^{(s)})|X_{ij} = x]$ and $H(x) =$ $E[(Z_{ij} - \mu_Z)|X_{ij} = x]$, respectively, and by $H_{a,mn}^{(s)}(x) \equiv \sum_{l=1}^{p} \hat{P}_{l,w}^{(s)}(x_k)$ and $H_{a,mn}(x) \equiv \sum_{l=1}^{p} \hat{P}_{l,w}^{Z}(x_l)$ the corresponding estimators of $H_a^{(s)}(x)$ and $H_a(x)$. Then we have

$$
\begin{aligned}
B_{mn}^{ZZ} = {} & \frac{1}{mn} \sum_{i=1}^{m} \sum_{j=1}^{n} \tilde{Z}_{ij}^{*} (\tilde{Z}_{ij}^{*})^{\tau} + \frac{1}{mn} \sum_{i=1}^{m} \sum_{j=1}^{n} \tilde{Z}_{ij}^{*} (\Delta_{ij}^{H_a})^{\tau} \\
& + \frac{1}{mn} \sum_{i=1}^{m} \sum_{j=1}^{n} \Delta_{ij}^{H_a} (\tilde{Z}_{ij}^{*})^{\tau} + \frac{1}{mn} \sum_{i=1}^{m} \sum_{j=1}^{n} \Delta_{ij}^{H_a} \Delta_{ij}^{H_a \tau} \\
\equiv {} & \sum_{k=1}^{4} B_{mn,k}^{ZZ},
\end{aligned}
\tag{A.25}
$$

where $\tilde{Z}_{ij}^{*} = \tilde{Z}_{ij} - H_a(X_{ij})$ and $\Delta_{ij}^{H_a} = H_a(X_{ij}) - H_{a,mn}(X_{ij})$. Moreover,

$$
\begin{aligned}
B_{mn}^{ZY} = {} & \frac{1}{mn} \sum_{i=1}^{m} \sum_{j=1}^{n} Z_{ij}^{*} \epsilon_{ij}^{*} + \frac{1}{mn} \sum_{i=1}^{m} \sum_{j=1}^{n} Z_{ij}^{*} (\Delta_{ij}^{(0)} - \Delta_{ij}^{H_a \tau} \beta) \\
& + \frac{1}{mn} \sum_{i=1}^{m} \sum_{j=1}^{n} \Delta_{ij}^{H_a} \epsilon_{ij}^{*} + \frac{1}{mn} \sum_{i=1}^{m} \sum_{j=1}^{n} \Delta_{ij}^{H_a} [\Delta_{ij}^{(0)} - (\Delta_{ij}^{H_a})^{\tau} \beta] \\
\equiv {} & \sum_{j=1}^{4} B_{mn,j}^{ZY},
\end{aligned}
\tag{A.26}
$$

where $\epsilon_{ij}^{*} = Y_{ij}^{*} - Z_{ij}^{* \tau} \beta$, $Z_{ij}^{*}$ and $Y_{ij}^{*} = \tilde{Y}_{ij} - H_a^{(0)}(X_{ij})$ are as defined in Assumption 3.2(i) and Theorem 3.1, and $\Delta_{ij}^{(s)} \equiv H_a^{(s)}(X_{ij}) - H_{a,mn}^{(s)}(X_{ij})$. So, to prove the asymptotic normality of $\hat{\beta}$, it suffices to show that

$$
B_{mn}^{ZZ} \xrightarrow{P} B^{ZZ}, \qquad \sqrt{mn}(B_{mn}^{ZY} - \mu_B) \xrightarrow{D} N(0, \Sigma_B),
\tag{A.27}
$$

where $B^{ZZ}$, $\mu_B$ and $\Sigma_B$ are as defined in Theorem 3.1. To this end, we need to have

$$
\sum_{i=1}^{m} \sum_{j=1}^{n} (\hat{P}_{k,w}^{(s)}(X_{ij}^{(k)}) - P_{k,w}^{(s)}(X_{ij}^{(k)}))^2 = o_P(\sqrt{mn}),
\tag{A.28}
$$

$$
s = 0, 1, \ldots, q.
$$

This is ensured by the following facts: due to (A.17), together with Lemma A.3 for $p = 1$,

$$
\sup_{x_k \in [-L_k, L_k]} |\hat{P}_{k,w}^{(s)}(x_k) - P_{k,w}^{(s)}(x_k)|
$$

$$= O_P((mnb_k^{1+2/r})^{-r/(1+2r)} + b_k^2) + O_P(1) \sum_{l=1,\neq k}^{p} b_l^2$$

$$+ O_P(1) \sum_{l=1,\neq k}^{p} b_l b_k + o_P(1) b_k^2 + O_P(1) \left( \frac{1}{\sqrt{mn}} \right),$$

and owing to $mnb_k^{4(2+r)/(2r-1)} \to \infty$ for some integer $r \geq 3$ and $mnb_k^5 = O(1)$,

$$\sqrt{mn}((mnb_k^{1+2/r})^{-r/(1+2r)} + b_k^2)^2$$

$$= C((mn)^{-(2r-1)/(1+2r)} b_k^{-4(2+r)/(1+2r)} + mnb_k^8)^{1/2}$$

$$\to 0,$$

$$\sqrt{mn} \left( O_P(1) \sum_{l=1,\neq k}^{p} b_l^2 + O_P(1) \sum_{l=1,\neq k}^{p} b_l b_k + o_P(1) b_k^2 + O_P(1) \left( \frac{1}{\sqrt{mn}} \right) \right)^2$$

$$\to 0.$$

Thus,

$$(A.29) \quad \sum_{i=1}^{m} \sum_{j=1}^{n} (\Delta_{ij}^{(s)})^2 = \sum_{i=1}^{m} \sum_{j=1}^{n} \left( \sum_{k=1}^{p} \hat{P}_{k,w}(X_{ij}^{(k)}) - P_{k,w}(X_{ij}^{(k)}) \right)^2$$

$$= o_P(\sqrt{mn}).$$

Therefore, using the Cauchy–Schwarz inequality, it follows that the $(s,t)$th element of $B_{mn,4}^{ZZ}$ satisfies

$$B_{mn,4}^{ZZ}(s,t) = \frac{1}{mn} \sum_{i=1}^{m} \sum_{j=1}^{n} \Delta_{ij}^{(s)} \Delta_{ij}^{(t)}$$

$$\leq \frac{1}{mn} \left( \sum_{i=1}^{m} \sum_{j=1}^{n} (\Delta_{ij}^{(s)})^2 \right)^{1/2} \left( \sum_{i=1}^{m} \sum_{j=1}^{n} (\Delta_{ij}^{(t)})^2 \right)^{1/2} = o_P(1),$$

and similarly

$$B_{mn,2}^{ZZ}(s,t) = o_P(1), \qquad B_{mn,3}^{ZZ}(s,t) = o_P(1).$$

Now since $B_{mn,1}^{ZZ} \to E[Z_{11}^* Z_{11}^{*\top}]$ in probability, it follows from (A.26) that the first limit of (A.27) holds with $B^{ZZ} = E[Z_{11}^* Z_{11}^{*\top}]$. To prove the asymptotic normality in (A.27), by using the Cauchy–Schwarz inequality and (A.29), we have

$$\sqrt{mn} \sum_{k=2}^{4} B_{mn,k}^{ZY} = o_P(1).$$

Therefore, the second limit of (A.27) follows from (A.26) and

$$\sqrt{mn}(B_{mn,1}^{ZY} - \mu_B) = \frac{1}{\sqrt{mn}} \sum_{i=1}^{m} \sum_{j=1}^{n} [Z_{ij}^* \epsilon_{ij}^* - \mu_B] \xrightarrow{D} N(0, \Sigma_B),$$

with $\mu_B = E[R_{ij}]$ and $\Sigma_B = \sum_{i=-\infty}^{\infty} \sum_{j=-\infty}^{\infty} E[R_{00} R_{ij}^\tau]$, where $R_{ij} = Z_{ij}^* \epsilon_{ij}^*$. The proof of the asymptotic normality follows directly from the central limit theorem for mixing random fields (see Theorem 6.1.1 of [20], e.g.). When (1.2) holds, the proof of the second half of Theorem 3.1 follows trivially. $\square$

PROOF OF COROLLARY 3.1. Its proof follows from that of Theorem 3.1. $\square$

PROOF OF THEOREM 3.2. Note that

$$\widehat{\widehat{P}}_{k,w}(x_k) = \widehat{P}_{k,w}^{(0)}(x_k) - \widehat{\beta}^\tau \widehat{P}_{k,w}^Z(x_k)$$

given in (2.12) and that $P_{k,w}(x_k) = P_{k,w}^{(0)}(x_k) - \beta^\tau P_{k,w}^Z(x_k)$. Then

$$\widehat{\widehat{P}}_{k,w}(x_k) - P_{k,w}(x_k)$$
$$= [\widehat{P}_{k,w}^{(0)}(x_k) - P_{k,w}^{(0)}(x_k) - \beta^\tau(\widehat{P}_{k,w}^Z(x_k) - P_{k,w}^Z(x_k))] - (\widehat{\beta} - \beta)^\tau \widehat{P}_{k,w}^Z(x_k)$$
$$= P_{mn,1}(x_k) + P_{mn,2}(x_k).$$

For any $c = (c_0, C_1^\tau)^\tau \in \mathbb{R}^{1+q}$ with $C_1 = (c_1, \ldots, c_q)^\tau \in \mathbb{R}^q$, we note that, for $x_k \in [-L_k, L_k]$,

$$\sum_{s=0}^{q} c_s P_{k,w}^{(s)}(x_k) = c_0 P_{k,w}^{(0)}(x_k) + C_1^\tau P_{k,w}^Z(x_k)$$

$$= E[g^{**}(X_{ij}^{(-k)}, x_k)] w_{(-k)}(X_{ij}^{(-k)}),$$

where $g^{**}(x) = E[Y_{ij}^{**} | X_{ij} = x]$ with $Y_{ij}^{**} = c_0(Y_{ij} - \mu_Y) + C_1^\tau(Z_{ij} - \mu_Z)$, and similarly,

$$\sum_{s=0}^{q} c_s \widehat{P}_{k,w}^{(s)}(x_k) = c_0 \widehat{P}_{k,w}^{(0)}(x_k) + C_1^\tau \widehat{P}_{k,w}^Z(x_k)$$

$$= \frac{1}{mn} \sum_{i=1}^{m} \sum_{j=1}^{n} g_{m,n}^{**}(X_{ij}^{(-k)}, x_k) w_{(-k)}(X_{ij}^{(-k)}),$$

where $g_{m,n}^{**}(x)$ is the local linear estimator of $g^{**}(x)$, as defined in Section 2 with $\widetilde{Y}_{ij}^{**} = c_0 \widetilde{Y}_{ij} + C_1^\tau \widetilde{Z}_{ij}$ instead of $\widetilde{Y}_{ij}$ there. Therefore, using the argument of Lemma A.5, the distribution of

$$\text{(A.30)} \qquad \sqrt{mnb_k} \sum_{s=0}^{q} c_s(\widehat{P}_{k,w}^{(s)}(x_k) - P_{k,w}^{(s)}(x_k))$$

is asymptotically normal.

Now taking $c_0 = 0$ in (A.30) shows that $\widehat{P}^Z_{k,w}(x_k) \to P^Z_{k,w}(x_k)$ in probability, which together with Theorem 3.1 leads to

$$(A.31) \quad \sqrt{mnb_k}P_{mn,2}(x_k) = \sqrt{mnb_k}(\widehat{\beta} - \beta)^\tau \widehat{P}^Z_{k,w}(x_k) = O_P(\sqrt{b_k}) = o_P(1).$$

On the other hand, taking $c_0 = 1$ and $C_1 = -\beta$ in (A.30), we have

$$(A.32) \quad \begin{aligned} &\sqrt{mnb_k}P_{mn,1}(x_k) \\ &\quad = \sqrt{mnb_k}[\widehat{P}^{(0)}_{k,w}(x_k) - P^{(0)}_{k,w}(x_k) - \beta^\tau(\widehat{P}^Z_{k,w}(x_k) - P^Z_{k,w}(x_k))] \end{aligned}$$

are asymptotically normal as in (A.8), with $Y^{**}_{ij} = Y_{ij} - \mu_Y - \beta^\tau(Z_{ij} - \mu_Z)$ and $g^{**}(x) = E(Y^{**}_{ij}|X_{ij} = x)$ instead of $H^{(s)}(x)$ and $Z^{(s)}_{ij}$ in Lemma A.5, respectively. This finally yields Theorem 3.2. $\square$

**Acknowledgments.** We would like to thank the Editors, an Associate Editor and the referees for their constructive comments and suggestions.

J. GAO
SCHOOL OF MATHEMATICS AND STATISTICS
UNIVERSITY OF WESTERN AUSTRALIA
CRAWLEY, WESTERN AUSTRALIA 6009
AUSTRALIA
E-MAIL: jiti@maths.uwa.edu.au

Z. LU
INSTITUTE OF SYSTEMS SCIENCE
ACADEMY OF MATHEMATICS
    AND SYSTEMS SCIENCE
CHINESE ACADEMY OF SCIENCES
BEIJING 100080
P. R. CHINA
AND
DEPARTMENT OF STATISTICS
LONDON SCHOOL OF ECONOMICS
HOUGHTON STREET
WC2A 2AE
UNITED KINGDOM
E-MAIL: z.lu@lse.ac.uk

D. TJØSTHEIM
DEPARTMENT OF MATHEMATICS
UNIVERSITY OF BERGEN
BERGEN 5007
NORWAY
E-MAIL: dagt@mi.uib.no